\documentclass[a4paper]{article}

\usepackage{amsmath,amssymb,amsthm,latexsym,graphicx,indentfirst,subfig, pb-diagram,wasysym,accents
}
\usepackage{fullpage}
\usepackage{calrsfs}
\usepackage{titling}
\title{Entropy numbers of diagonal operators on Orlicz sequence spaces}

\author{\large{Thanatkrit Kaewtem and Yuri Netrusov}
}

\date{}

\theoremstyle{theorem}
\newtheorem{thm}{Theorem}[section]
\newtheorem{lem}[thm]{Lemma}
\newtheorem{cor}[thm]{Corollary}
\newtheorem{prop}[thm]{Proposition}

\theoremstyle{definition}
\newtheorem{defn}[thm]{Definition}

\theoremstyle{remark} 
\newtheorem{rem}[thm]{Remark}

\newcommand{\R}{\mathbb{R}}
\newcommand{\N}{\mathbb{N}}
\newcommand{\Z}{\mathbb{Z}}

\begin{document}

\maketitle
	

\begin{abstract}
Let  $M_1$ and $M_2$ be functions on $[0,1]$ such that $M_1(t^{1/p})$ and $M_2(t^{1/p})$ are Orlicz functions for some $p \in (0,1].$ 
Assume that 
$M_2^{-1} (1/t)/M_1^{-1} (1/t)$ is non-decreasing for $t \geq 1.$
Let $(\alpha_i)_{i=1}^\infty$ be a non-increasing sequence of non-negative real numbers. 
Under some conditions on $(\alpha_i)_{i=1}^\infty,$ sharp two-sided estimates for entropy numbers of diagonal operators $T_\alpha :\ell_{M_1} \rightarrow \ell_{M_2}$ generated by $(\alpha_i)_{i=1}^\infty,$ where $\ell_{M_1}$ and $\ell_{M_2}$ are Orlicz sequence spaces, are proved. 
The results generalise some works of Edmunds and Netrusov in \cite{E&N2010} and hence a result of Cobos, K\"{u}hn and Schonbek in \cite{Cobos&Kuhn&Schonbek}.
\end{abstract}

\section{Introduction}

Let $X, Y$ be real quasi-Banach spaces and $T$ a continuous linear operator from $X$ into $Y.$
For $k \in \N,$ the  $k^{\text{th}}$ (dyadic) {\it entropy number} of $T$ is defined by
\[ e_k(T) = \inf \left\lbrace \varepsilon > 0 : T(B_X) \subset \bigcup_{i=1}^{2^{k-1}}(y_i + \varepsilon B_Y) ~\text{for some}~ y_1, y_2, ..., y_{2^{k-1}} \in Y \right\rbrace,\]
where $B_X$ and $B_Y$ denote the closed unit balls in $X$ and $Y,$ respectively.	 
It can be seen that the numbers $e_k (T)$ are monotonic decreasing as $k$ increases, with $e_1(T) \leq \|T\|.$
In addition, $T$ is compact if and only if $\displaystyle \lim_{k \rightarrow \infty} e_k(T) = 0.$
Hence, the rate of decrease may be regarded as a measure of the degree of compactness of $T.$
For the background and applications of entropy numbers, we refer to any textbooks given in the references, especially  \cite{C&S, E&T,P0}.


Entropy numbers of embeddings between function spaces have been of interest over the past decades.
A pioneer result in the area, due to Kolmogorov and Tikhomirov in 1959, is the study of entropy numbers of embeddings of $C^{m}([0,1]^n)$ in $C([0,1]^n),$ where $m\in \N$ (see \cite{Kol&Ti}).
In 1967, Birman and Solomyak proved a remarkable theorem concerning entropy numbers of Sobolev embeddings (see \cite{Birman&Solomyak1}).
Generalisations of these results can be found in the book written by Edmunds and Triebel \cite{E&T} and the references given therein.
Later, Haroske, Triebel, K\"{u}hn, Leopold, Sickel and Skrzypczak studied  the problem of estimating entropy numbers of embeddings of weighted function spaces of Besov and Sobolev type (see \cite{Haroske&Triebel1994}, \cite{Ku&Leo&Sick&Skr1, Ku&Leo&Sick&Skr2, Ku&Leo&Sick&Skr3}).
A technique which is usually used to estimate entropy numbers of function space embeddings is reduction to a simpler problem in sequence spaces.

Let $(\alpha_i)_{i=1}^\infty$ be a sequence of non-negative real numbers. 
A {\it diagonal operator} generated by $(\alpha_i)_{i=1}^\infty$ is the continuous linear operator $T_{\alpha}:X \rightarrow Y$ defined by, for each $i \in \N,$ $T_{\alpha}(x_i) = \alpha_i y_i,$  where $\{x_i \}_{i=1}^\infty$ and $\{y_i\}_{i=1}^\infty$ are bases for $X$ and $Y,$ respectively.
If $X$ is a Banach space with a 1-unconditional basis $\{x_i \}_{i=1}^\infty,$  then sharp estimates for the entropy numbers of the diagonal operator $T_{\alpha}:X \rightarrow X$ generated by a non-increasing sequence of non-negative real numbers $(\alpha_i)_{i=1}^\infty$  were given by Gordon, K\"{o}nig and Sch\"{u}tt in 1987 (see \cite{Gor&Ko&Sch}, p.226). 
The exact asymptotic behaviour of entropy numbers of general diagonal operators is still unknown when $X \neq Y.$  


Let $k \in \N, ~0 < p < q \leq \infty$ and $T_{\alpha}: \ell_{p} \rightarrow \ell_{q}$ a diagonal operator geneated by a non-increasing sequence of non-negative real numbers $(\alpha_i)_{i=1}^\infty.$
One can see that if $\alpha_1 = \alpha_2 =... = \alpha_n = 1$ and $\alpha_i = 0$ for $i > n$, then $T_{\alpha}$ is the projection $P_n$ onto the first $n$ coordinates.
In order to estimate $e_k(T_{\alpha}),$ the above special case suggests us to represent the map as $T_{\alpha} = T_{\sigma} + T_{\mu},$ where
$T_{\sigma}$ and $T_{\mu}$ are diagonal operators generated by $(\sigma_{i})_{i=1}^\infty$ and $(\mu_{i})_{i=1}^\infty$  defined by $\sigma_i = \min \{ \alpha_i, \alpha_k \}$ and $\mu_i = \alpha_i - \sigma_i,$ respectively, and then to estimate the entropy numbers of these two operators.	
In 2010, Edmunds and Netrusov established sharp two-sided estimates of $e_k(T_{\sigma})$ (see \cite{E&N2010}).
The result gave directly precise estimates for the entropy numbers $e_k(P_n)$ for all $k \leq n$, or, equivalently, for the entropy
numbers of identity maps between finite-dimensional Lebesgue sequence spaces which were discovered by Sch\"{u}tt (see \cite{Sch}) in the Banach space setting.
The sharp two-sided estimates of $e_k(T_{\mu})$ are still unknown, however, under mild regularity and decay conditions on the generating sequence $(\alpha_i)_{i=1}^\infty$, the entropy numbers $e_k(T_{\mu})$ can be estimated, and hence, in combination with knowledge of $e_k(T_{\sigma})$,  the entropy numbers of the given diagonal operator $e_k(T_\alpha)$ can also be obtained.
This will be illustrated in Corollary \ref{Cor-Doubling-Orlicz-sp-Q} in our paper.

Generalisations of $\ell_p$-spaces were introduced by many mathematicians, for example, G. Lorentz and W. Orlicz  (as cited in the monograph \cite{Lind&Tzaf}). 
Orlicz sequence spaces were inspired by the function $t^p$ in the definition of $\ell_p$-spaces or, more generally $L_p$.   
It is natural to try to replace $t^p$ by a more general function $M,$ and to consider the set of all sequences $(x_i)_{i=1}^\infty$ for which the series $\sum_{i=1}^{\infty} M(|x_i|)$ converges. 
There are some advantages of Orlicz spaces but we are not going to discuss these in this paper.
Orlicz (see \cite{Orlicz}) provided the restrictions on the function $M$ in order to make this set of sequences into a Banach space. It is well-known that Orlicz sequence spaces can be extended to the case of quasi-Banach spaces (see, for instance, Remark  \ref{Rmk-l_M-is-Quasi-B-Sp} and Proposition \ref{Prop-Gamma-norm-on-Orlicz} below).

While the entropy behaviour of some diagonal maps acting between $\ell_p$-spaces is known, it was the lack of similar knowledge when the sequence spaces are more sophisticated that was responsible for a delay in establishing such matters as the poor behaviour of entropy numbers under real interpolation (see \cite{E&N2010}).  
This paper furthers our understanding in the natural case when the spaces are Orlicz sequence spaces.

Let  $M_1$ and $M_2$ be functions on $[0,1]$ such that $M_1(t^{1/p})$ and $M_2(t^{1/p})$ are Orlicz functions for some $p \in (0,1].$ 
Then the sequence spaces $\ell_{M_1}$ and $\ell_{M_2}$ are $p$-Banach spaces (see Proposition \ref{Prop-Gamma-norm-on-Orlicz}).
Denote by $\varphi_{M_1}$ and $\varphi_{M_2}$ the fundamental functions of $\ell_{M_1}$ and $\ell_{M_2},$ respectively (see Definition \ref{Defi-Fundamental-fn}).
Let $(\alpha_i)_{i=1}^\infty$ be a non-increasing sequence of non-negative real numbers and $T_\alpha :\ell_{M_1} \rightarrow \ell_{M_2}$ the diagonal operator generated by $(\alpha_i)_{i=1}^\infty$.
In this paper, using combinatorial techniques, we prove that (see Theorem \ref{GeneralEdm&Yu-Orlicz-sp-Q} and Proposition \ref{Rem-Equiv-of-Theta-Q}) if $\varphi_{M_1}/\varphi_{M_2}$ is a non-decreasing function and $(\alpha_i)_{i=1}^\infty$ satisfies  the condition that $\alpha_1 = \alpha_2 =...=\alpha_{k}$ for some $k \in \N,$ then 
\begin{equation} \label{Eq1-Introduction}
	c_1(p) \max_{s \in \{1,...,k\}} \alpha_{k2^{s-1}} \frac{\varphi_{M_2} \left( k/s\right) }{\varphi_{M_1} \left( k/s\right) }
	\leq	e_k(T_\alpha ) 
	\leq  c_2(p) \max_{s \in \{1,...,k\}} \alpha_{k2^{s-1}} \frac{\varphi_{M_2} \left(k/s\right) }{\varphi_{M_1} \left( k/s\right) } 
\end{equation}	 
where $c_1(p) = 4^{-2}  2^{-21/p }  \left(  \log \left( \frac{24}{p} \right) \right)^{-2/p} \left( \frac{1}{p} \right) ^{-1/p }$ and 
$c_2(p) = 4 \cdot 2^{20/p }  \left(  \log \left( \frac{24}{p} \right) \right)^{2/p} \left( \frac{1}{p} \right) ^{1/p }.$
The lower estimate of (\ref{Eq1-Introduction}) follows easily from the extension of Sch\"{u}tt's result (see \cite{Sch}, Theorem 1) when the natural embedding is defined on Orlicz sequence spaces (see Theorem \ref{Thm-Schutt'sThmForOrliczSp-New1} and also Theorem \ref{Thm-Schutt'sThmForOrliczSp-New2}). The upper estimate of (\ref{Eq1-Introduction}), on the other hand, is an application of a generalisation of Lemma 11 proved by Edmunds and Netrusov in \cite{E&N2010} (see Lemma \ref{Lemma-Before-main-Orlicz-seq-spaces-Q}).
The inequalities (\ref{Eq1-Introduction}) is an extension of Theorem 12 in \cite{E&N2010} where no computation of constants was given. 
This is a new result even in the Banach space case ($p = 1 $).
We note that the constants $c_1(p)$ and $c_2(p)$ above can be improved.
In the special case when $T_\alpha$ is an embedding between finite dimensional Lebesgue sequence spaces, the estimates of constants were given by Gu\'{e}don and Litvak in \cite{Gued&Litv}.
Finally, we apply the inequalities (\ref{Eq1-Introduction}) to give sharp two-sided estimates for the entropy numbers of diagonal operators $T_\alpha :\ell_{M_1} \rightarrow \ell_{M_2}$ generated by a non-negative non-increasing sequence $(\alpha_i)_{i=1}^\infty$ which satisfies the doubling condition (see Corollary \ref{Cor-Doubling-Orlicz-sp-Q}).
\section{Definitions and notations}

In this section we collect some basic facts, conventions and definitions that will be used later. 
All spaces considered will be assumed to be real vector spaces.
For $p \in (0,1],$ a {\it $p$-norm} on a vector space $X$ is a function $\|\cdot \|_{X}:X \rightarrow [0,\infty)$ which has the property of a norm, but instead of the triangle inequality, the inequality 
\begin{equation} \label{Eq-DefofGammaNorm}
	\|x+y \|_{X}^p \leq \|x \|_{X}^p + \|y\|_{X}^p
\end{equation}
holds for any $x,y \in X$.
It is clear that every $q$-norm is a $p$-norm for $0<p \leq q\leq 1.$
If (\ref{Eq-DefofGammaNorm}) is replaced by 
\begin{equation}  \label{Eq-DefofQuasiNorm}
	\|x+y\|_{X} \leq C(\|x \|_{X}+ \|y \|_{X})
\end{equation}
for some constant $C \geq 1$, then it is called a {\it quasi-norm}. 
It is also clear that every $p$-norm is a quasi-norm with $C = 2^{1/p-1}.$ 
According to the Aoki-Rolewicz theorem (\cite{Aoki} and \cite{Rolewicz}), if $X$ is a quasi-Banach space with a constant $C$, then there exists $p \in (0,1]$ (in fact, $ 2^{1/p-1} = C$) and a $p$-norm on $X$ which is equivalent to the original quasi-norm. 

Given quasi-Banach spaces $X$ and $Y$, we write $\mathcal{B}(X,Y)$ for the space of all continuous linear maps from $X$ into $Y,$ and write $\mathcal{B}(X)$ if $X=Y.$ 
Logarithms are always taken to be base 2, $\log = \log_2,$ and $\lfloor x \rfloor$ will denote the integer part of the real number $x.$
We will write $c:=c(a_1,a_2,...,a_n)$ to emphasise that $c$ depends only on the parameters $a_1,a_2,...,a_n.$

The following properties of entropy numbers will be used later. 
If $Y$ is a $p$-Banach space, then, for each $k_1,k_2 \in \N,$
\[e_{k_1+k_2-1}^p (T_1+T_2) \leq e_{k_1}^p (T_1) + e_{k_2}^p (T_2) ~\text{and}~ e_{k_1+k_2-1} (RS) \leq e_{k_1}(R) e_{k_2}(S),\]
whenever $T_1+T_2$ and $RS$ are properly defined operators (see \cite{E&T}).
The  $k^{\text{th}}$ (dyadic) {\it inner entropy number} $f_k(T)$ of $T \in \mathcal{B}(X,Y)$ is defined to be the supremum of all those $\varepsilon > 0$ such that there are $x_1,...,x_{2^{k-1}+1} \in B_X$ with $\|Tx_i - Tx_j \|_{Y} \geq 2 \varepsilon$ whenever $i,j$ are distinct elements of $\{1,...,2^{k-1}+1\}.$
If $Y$ is a $p$-Banach space, then the outer and inner entropy numbers are related by
\begin{equation} \label{EqRelationInner&Outer}
	2 ^{1-1/p}f_k(T) \leq e_k(T) \leq 2f_k(T).
\end{equation}
These estimates were proved by Pietsch in the Banach space case ($p = 1$); a simple modification gives us (\ref{EqRelationInner&Outer}). 


\begin{defn} \label{Defi-Orlicz-fn}
	An {\it Orlicz function} $N$ on $[0,1]$ is a continuous, strictly increasing and convex function defined on $[0,1]$ such that $N(0) = 0$ and $N(1) = 1.$ 
\end{defn}


\begin{defn}
	Let $M$ be a function on $[0,1]$ such that $M(t^{1/p})$ is an Orlicz function for some $p \in (0,1].$   
	The space $\ell_{M} := \{x = (x_i)_{i=1}^{\infty} :  \sum_{i=1}^{\infty} M(|x_i|/\varrho) < \infty  ~\text{ for some}~ \varrho > 0 \}$ equipped with the Luxemburg norm 
	\[ \|x\|_{\ell_M} := \inf \left\lbrace  \varrho >0 : \sum_{i=1}^{\infty} M \left( \frac{|x_i|}{\varrho} \right)  \leq 1 \right\rbrace \]
	is a $p$-Banach space (see Proposition \ref{Prop-Gamma-norm-on-Orlicz} below) called an {\it Orlicz sequence space}. 
\end{defn}
\noindent Clearly, if $M(t) = t^p$ for some $0< p<\infty$, then $\ell_M$ is just the sequence space $\ell_p$ with its usual norm.

The next simple result can be proved analogously as in the Banach space case.

\begin{prop} \label{Prop-Gamma-norm-on-Orlicz}
	Let $M$ be a function on $[0,1]$ such that $M(t^{1/p})$ is an Orlicz function for some $p \in (0,1].$  
	Then $\|\cdot\|_{\ell_{M}}$ is a $p$-norm on $\ell_{M}.$
\end{prop}
\begin{proof}
	Let $N$ be the convex function on $[0,1]$ defined by $N(t) = M(t^{1/p}).$
	Let $x = (x_i)_{i=1}^\infty, y = (y_i)_{i=1}^\infty \in \ell_{M}.$ 
	Set $\sigma_1 = \| x\|_{\ell_M}, \sigma_2 = \| y\|_{\ell_M}$ and $\varrho = (\sigma_1^p + \sigma_2^p)^{1/p}.$
	Without loss of generality, we may assume that $\sigma_1,\sigma_2 \neq 0.$
	Then,
	\begin{align*}
		\sum_{i=1}^{\infty} M\left( \frac{|x_i +y_i|}{\varrho}\right)
		&= \sum_{i=1}^{\infty} N\left( \frac{|x_i +y_i|^p}{\varrho^p}\right) \\
		&\leq \sum_{i=1}^{\infty}  N\left( \frac{|x_i |^p}{\sigma_1^p} \cdot \frac{\sigma_1^p}{\varrho^p} + \frac{|y_i |^p}{\sigma_2^p} \cdot \frac{\sigma_2^p}{\varrho^p} \right) \\
		&  \leq\frac{\sigma_1^p}{\varrho^p} \sum_{i=1}^{\infty}  N\left( \frac{|x_i |^p}{\sigma_1^p}\right) + \frac{\sigma_2^p}{\varrho^p} \sum_{i=1}^{\infty}  N\left( \frac{|y_i |^p}{\sigma_2^p}\right) \\
		&  =   \frac{\sigma_1^p}{\varrho^p} \sum_{i=1}^{\infty}  M\left( \frac{|x_i |}{\sigma_1}\right) + \frac{\sigma_2^p}{\varrho^p} \sum_{i=1}^{\infty}  M\left( \frac{|y_i |}{\sigma_2}\right)  \\
		&\leq 1.
	\end{align*}
	This implies that $\|x+y\|_{\ell_M} \leq  \varrho =  \left( \| x\|_{\ell_M}^p + \| y\|_{\ell_M}^p\right) ^{1/p}.$
\end{proof}

Let $E$ be a non-empty finite or countable set.
Given any function $x:E \rightarrow \R$, we write
\[\|x \|_{\ell_{M}(E)} = \inf \left\lbrace  \varrho >0 : \sum_{j \in E} M \left( \frac{|x(j)|}{\varrho} \right)  \leq 1 \right\rbrace,\]
and denote by $\ell_{M}(E) := \{x:E \rightarrow \R : \sum_{j \in E} M(|x(j)|/\varrho) < \infty  ~\text{ for some}~ \varrho > 0\}$.  
We will write $\|\cdot\|_{\ell_{M}} $ instead of $\|\cdot\|_{\ell_{M}(E)} $ if no ambiguity is possible.
The space $\ell_{M}(E)$ is a quasi-Banach space. 
The elements of the space $\ell_{M}(E)$ may be denoted by the sequence $\{x(j)\}_{j \in E}$ of their values on $E$, and are often written as $(x_j)_{j \in E}$. 
In particular, when $E=\N$ or $E = \{1,2,...,n\}$ for some $n \in \N$, we will denote $\ell_{M}(E)$ by $\ell_{M}$ or $\ell_{M}^n,$ respectively.

\begin{defn} \label{Defi-Fundamental-fn}
	Let $M$ be a function on $[0,1]$ such that $M(t^{1/p})$ is an Orlicz function for some $p \in (0,1].$   
	A non-negative function $\varphi_{M}$ defined on $[1,\infty)$ by
	\begin{equation}
		\varphi_{M} (t)  = \frac{1}{M^{-1} (1/t)}
	\end{equation}
	is called a {\it fundamental function} of $\ell_{M}.$
	Here, $M^{-1}$ stands for the inverse function of $M.$
\end{defn}
\begin{rem}
	Let $M$ be a function on $[0,1]$ such that $M(t^{1/p})$ is an Orlicz function for some $p \in (0,1].$   
	\begin{enumerate}
		\item		For any $n \in \N,$ it can be seen that
		\begin{equation}
			\varphi_M(n) = \|e_1+e_2+...+e_n\|_{\ell_M},
		\end{equation}
		where $e_1,e_2,...,e_n$ are standard unit vectors in $\ell_M.$
		\item The fact that the map $t \mapsto \frac{M^{-1}(t)}{t^{1/p}}$ is decreasing for $t \in (0,1]$ is equivalent to the statement that the map $t \mapsto \frac{t^{1/p}}{\varphi_{M}(t)}$ is increasing for $t \in [1,\infty).$  
		Hence,
		\[\varphi_{M}(t) \leq \varphi_{M}(Ct) \leq C^{1/p} \varphi_{M}(t).\]
		for any $C \geq 1$ and $t \geq 1.$
	\end{enumerate}
\end{rem}



From now on, we let  $M_1$ and $M_2$ be functions on $[0,1]$ such that $M_1(t^{1/p})$ and $M_2(t^{1/p})$ are Orlicz functions for some $p \in (0,1].$ In addition, we denote by $\varphi_{M_1}$ and $\varphi_{M_2}$ the fundamental functions of $\ell_{M_1}$ and $\ell_{M_2},$ respectively.


\section{Preliminary results}


This section presents some useful results that will be needed.
Many results have a quite technical appearance due to the attempts at computing constants.

First, we give two-sided estimates for entropy numbers of the identity map between finite-dimensional $p$-Banach spaces. 
The proof is similar to the Banach case (see, for instance, \cite{C&S}). 
\begin{thm}  \label{Entorpy-of-finite-dim-identity-map-Q}
	Let $k \in \N$ and $X$ an $n$-dimensional $p$-Banach space. If $I:X \rightarrow X$ is the identity map, then 
	\begin{equation*}
		2^{\frac{1-k}{n}} \leq e_k(I:X \rightarrow X) \leq 4^{1/p}\cdot 2^{\frac{1-k}{n}}.
	\end{equation*}
\end{thm}

The following well-known result can be found, for example, in \cite{Abra&Steg} page 257.
\begin{thm}[Stirling's formula]
	For any $n \in \N,$
	\begin{equation} \label{Stirlings-formula}
		\sqrt{2\pi} n^{n+1/2} e^{-n} \leq n! \leq \sqrt{2\pi} n^{n+1/2} e^{-n} e^{1/12n}.
	\end{equation}
\end{thm}
Hence, for any $n \in \N$ and $\beta \in (0,1),$ if $\beta n \in \N,$ by Stirling's formula,  the following estimate holds:
\begin{equation}
	{n \choose \beta n} \leq \frac{1}{\sqrt{2 \pi \beta (1-\beta)}} \frac{e^{1/12n}}{\sqrt{n}} \left[ \left(\frac{1}{\beta} \right)^\beta \left(\frac{1}{1-\beta} \right)^{1-\beta}  \right]^n. 
\end{equation}

\begin{lem} \label{Lem-Comnatorial-5}
	Let $m$ be a positive integer such that $m \geq 5$ and $n = 2^m.$
	Then there exists a set $\Omega \subseteq \Z^n$ which has the following properties:
	\begin{enumerate}
		\item[\textnormal{(i)}] $\sharp \Omega \leq 2^{\frac{n}{2}-1};$
		\item[\textnormal{(ii)}]  for any sequence $(x_i)_{i=1}^{n} \in \R^n$  such that
		\begin{equation}
			\forall j \in \{0,1,...,m-5\}, ~ x_{\lfloor 2^{j-1} \rfloor +1}^{\ast} \leq 2^{\psi(m,j)}  ,						
		\end{equation}
		where $\psi(m,j) = 2^{m-j-4} (m-j)^{-2},$
		there is a sequence $(z_i)_{i=1}^n \in \Omega$  such that 
		\begin{equation*}
			\max_{1 \leq i \leq n}|x_i-z_i| \leq 4 .
		\end{equation*}	
		Here the sequence $(x_i^\ast)_{i=1}^n$ is the non-increasing rearrangement of $(|x_i|)_{i=1}^n.$
	\end{enumerate}
\end{lem}

\begin{proof}
	Let $Q$ be the set of all sequences $(z_i)_{i=1}^{n/32} \in \Z^{n/32}$ such that  $z_1,...,z_{n/32}$ are divisible by 4  and 
	\begin{equation*}
		\forall j \in \{0,1,...,m-5\} \forall i \in \N 
		\left( 2^{j-1} < i \leq 2^{j} ~\Longrightarrow~   |z_i| \leq 2^{ \psi(m,j) } \right) .
	\end{equation*}
	We note that $1+2 \lfloor k/4 \rfloor \leq k$ for all $k \in \N.$
	If $m=5$ then
	\begin{equation} \label{Eq0-Lem-Comnatorial-5}
		\sharp Q \leq 2^{\psi(5,0)} = 2^{2/25} \leq 2^{1/4}.
	\end{equation}
	Now if $m \geq 6,$ it follows from the definition of the set $Q$ that
	\begin{align} \label{Eq1-Lem-Comnatorial-5}
		\log \sharp Q 
		&\leq \log  \left(  2^{\psi(m,0)} \right) 
		\left( 2^{\psi(m,1)} \right) 
		\left(  2^{\psi(m,2)}\right)^{2^{1}} 
		\left(  2^{\psi(m,3)}\right)^{2^{2}} \cdots 
		\left( 2^{\psi(m,m-5)} \right)^{2^{m-6}}  \notag \\
		&= 2^{m-5} \left( \frac{2}{m^2} + \frac{1}{(m-1)^2} +\frac{1}{(m-2)^2} +...+\frac{1}{5^2} \right)  \notag \\
		&\leq 2^{m-5} \left( \frac{1}{m^2} + \frac{1}{m-1} - \frac{1}{m} +\frac{1}{m-2}-\frac{1}{m-1} +...+\frac{1}{4}-\frac{1}{5} \right)   \notag \\
		&\leq 2^{m-5} \left( \frac{1}{m^2} -\frac{1}{m} + \frac{1}{4} \right)  \notag\\
		&\leq 2^{m-7}.
	\end{align}	
	Let $\Gamma$ be a set given by 
	\begin{equation*}
		\Gamma = \left\lbrace (E_j)_{j=0}^{m-5} : \sharp E_j = 2^j ~\text{for}~ j = 0,1,...,m-5 ~\text{and}~ E_0 \subset E_1 \subset...\subset E_{m-5} \subset \{1,...,n\} \right\rbrace.
	\end{equation*}
	Applying Stirling's formula to the first term in the product, we obtain the following estimates
	\begin{align} \label{Eq2-Lem-Comnatorial-5}
		\sharp \Gamma 
		&\leq {2^{m} \choose 2^{m-5}} {2^{m-5} \choose 2^{m-6}} {2^{m-6} \choose 2^{m-7}} \cdots {2^{2} \choose 2^{1}}  {2^{1} \choose 2^{0}}  \notag \\
		&\leq \frac{1}{\sqrt{2 \pi}} \frac{2^5}{\sqrt{2^{5}-1}} \frac{e^{1/384}}{2^{5/2}} \left( 32^{1/32} \left( 1+\frac{1}{31} \right)^{31/32} \right)^{2^m}  
		2^{2^{m-5}}  2^{2^{m-6}} \cdots  2^{2^{2}}  2^{2^{1}} \notag \\
		&\leq 
		\frac{1}{2} \left(  2^{5/32} e^{1/32} \right) ^{2^m}    2^{2^{m-4}}.
	\end{align}
	Now let $\Omega$ be the set of all sequences $(z_i)_{i=1}^{n} \in \Z^{n}$ such that  $z_1,...,z_{n}$ are divisible by 4  and  there exists a sequence $(E_j)_{j=0}^{m-5} \in \Gamma$ such that
	\[
	\begin{cases}
		|z_i| \leq 2^{\psi(m,0)} &~~\text{if}~~ i=1\\
		|z_i| \leq 2^{\psi(m,j)} &~~\text{if}~~ i \in E_j \setminus E_{j-1}, ~ j = 1,2,...,m-5\\
		z_i = 0 &~~\text{if}~~ i \in \{1,2,...,n\} \setminus E_{m-5}.
	\end{cases}
	\]
	It follows from (\ref{Eq0-Lem-Comnatorial-5}), (\ref{Eq1-Lem-Comnatorial-5}) and (\ref{Eq2-Lem-Comnatorial-5}) that 
	\[\sharp \Omega 
	\leq 
	\frac{1}{2} \left(  2^{\frac{1}{128} +\frac{1}{16} + \frac{5}{32}} e^{\frac{1}{32}}    \right)^{2^m} 
	\leq 2^{(2^{m-1} )-1}.
	\]
	Therefore, the property (i) holds.
	The property (ii) follows from the construction of $\Omega.$
\end{proof}
The following estimate will be used to prove the next theorem.

\begin{lem} \label{Lem-Estimate-using-Calulus}
	Let $0<\beta \leq 1$ and and set $A_{\beta} := \beta^{-1}16\log e.$
	Define the function $f:[2^5, \infty) \rightarrow (0,\infty)$  by 
	\[f(t) = t^{1/\beta} 2^{-\frac{t}{16(\log t)^2}}.\]
	Then 	
	$f(t) \leq 
	\left( 22A_{\beta}(\log A_{\beta})^2\right) ^{1/\beta}
	2^{-\frac{A_{\beta}(\log A_{\beta})^2}{16\left( \log \left( A_{\beta}(\log A_{\beta})^2 \right) \right) ^2}} 
	$
	for all $t \in [2^5, \infty).$
\end{lem}

\begin{proof}
	It is easy to verify that
	\[ f^{\prime}(t) = t^{1/\beta} 2^{-\frac{t}{16(\log t)^2}} \left[ \frac{(16 \log e) (\log t)^3 - \beta t \left( \log t - 2\log e\right)  }{\beta (16 \log e)t  (\log t)^3} \right]. \]
	Moreover,
	(i) $f^{\prime}(t) \geq 0$ if $2^5 \leq t \leq A_{\beta}(\log A_{\beta})^2$ and
	(ii) $f^{\prime}(t) \leq 0$ if $t \geq 22 A_{\beta}(\log A_{\beta})^2.$
	The desired result follows from (i) and (ii).
\end{proof}

\begin{thm} \label{Thm-Polynomial-decay-Orlicz-F}
	Let $m \in \N, m \geq 5$ and $\alpha \geq 0.$
	Let $D_{\alpha} : \ell_{M_1}^{2^m} \rightarrow \ell_{M_2}^{2^m}$ be the linear map defined by 
	\[ D_{\alpha}(x_i)_{i=1}^{2^m}= \left( \left( \frac{2^m}{i} \right)^{\alpha} x_i\right) _{i=1}^{2^m}.\]
	Then,
	\begin{equation} \label{Eq0-Polynomial-decay-Orlicz-F}
		e_{2^{m-1}}\left( D_{\alpha}: \ell_{M_1}^{2^m} \rightarrow \ell_{M_2}^{2^m}\right)  \leq c(\alpha, p)  \frac{\varphi_{M_2}(2^m)}{\varphi_{M_1}(2^m)},
	\end{equation}
	where $c(\alpha, p) = 4 \cdot 2^{11(1/p +\alpha)} \left(  \log \left( 24(1/p +\alpha)\right) \right)^{2(1/p +\alpha)} \left( 1/p +\alpha \right) ^{1/p +\alpha}.$
\end{thm}

\begin{proof}
	For any $n \in \N$ and $0<p<\infty,$ we denote by $\ell_{p,\infty}^n$ the $n$-dimensional Lorentz sequence space equipped with the quasi-norm
	$\displaystyle \|(x_i)_{i=1}^{n}\|_{\ell_{p,\infty}^n} = \sup_{1 \leq i \leq n} i^{1/p} x_i^\ast,$ where  $(x_i^\ast)_{i=1}^n$ is the non-increasing rearrangement of $(|x_i|)_{i=1}^n.$
	First, we are going to show that 
	\begin{equation} \label{Eq1-Polynomial-decay-Orlicz-F}
		D_{\alpha} \left( \varphi_{M_1}(2^m) B_{\ell_{M_1}^{2^m}} \right) 
		\subseteq 2^{1/\beta} 2^{m/\beta} B_{\ell_{\beta, \infty}^{2^m}}, 
	\end{equation}
	where $1/\beta = \alpha + 1/p.$
	Let $x = (x_i)_{i=1}^{2^m} \in B_{\ell_{M_1}^{2^m}}.$
	For any fixed $j \in \{1,2,3,...,2^m\},$ we have that $\varphi_{M_1} (2^m) 
	\leq \left( \frac{2^m}{j}\right)^{1/p} \varphi_{M_1} (j).$ 
	In addition, since $\|x\|_{\ell_{M_1}^{2^m}} \leq 1,$ it turns out that
	$j M_1 (x_j^\ast) \leq \sum_{i=1}^{j} M_1 (x_i^\ast)  
	\leq 1,$ and hence 
	\[ j^{1/p}  x_j^\ast \leq j^{1/p}  M_1^{-1}(1/j) = \frac{j^{1/p} }{\varphi_{M_1}(j)} \leq  \frac{2^{m/p}}{\varphi_{M_1}(2^m)}. \]
	This implies that $\|x\|_{\ell_{p,\infty}^{2^m}} \leq  \frac{2^{m/p}}{\varphi_{M_1}(2^m)}.$
	Consequently,
	\begin{equation} \label{Eq2-Polynomial-decay-Orlicz-F}
		D_{\alpha} \left( \varphi_{M_1}(2^m) B_{\ell_{M_1}^{2^m}} \right) 
		\subseteq D_{\alpha} \left( 2^{m/p} B_{\ell_{p,\infty}^{2^m}} \right). 
	\end{equation}
	Now the inclusion (\ref{Eq1-Polynomial-decay-Orlicz-F}) will be proved if we show that  $D_{\alpha} \left( 2^{m/p} B_{\ell_{p,\infty}^{2^m}} \right) \subseteq 2^{1/\beta} 2^{m/\beta} B_{\ell_{\beta, \infty}^{2^m}}.$ 
	To do this, let us fix $z = (z_i)_{i=1}^{2^m} \in B_{\ell_{p,\infty}^{2^m}}.$ 
	For any $(u_i)_{i=1}^{2^m}, (v_i)_{i=1}^{2^m} \in \ell_{\beta,\infty}^{2^m},$ we notice that  if $1 \leq k \leq 2^{m-1}$ then
	$(u_{2k-1} v_{2k-1})^{\ast} \leq u^{\ast}_{k} v^{\ast}_{k}$  (see, for example, \cite{P1} page 74).
	It follows that
	\begin{align*}
		\| (i^{-\alpha} z_i)_{i=1}^{2^m} \|_{\ell_{\beta,\infty}^{2^m}}
		&=\max  \left\lbrace \sup_{1 \leq k \leq 2^{m-1}} (2k-1)^{1/\beta} \left( \frac{z_{2k-1}}{(2k-1)^{\alpha}} \right)^{\ast}, \sup_{1 \leq k \leq 2^{m-1}} (2k)^{1/\beta} \left( \frac{z_{2k}}{(2k)^{\alpha}} \right)^{\ast}  \right\rbrace \\
		&\leq  \sup_{1 \leq k \leq 2^{m-1}} (2k)^{1/\beta} \left( \frac{z_{2k-1}}{(2k-1)^{\alpha}} \right)^{\ast}  \\
		&\leq 2^{1/\beta} \sup_{1 \leq k \leq 2^{m-1}} k^{1/\beta-\alpha}  z_{k}^{\ast} \\
		&\leq 2^{1/\beta} \| z\|_{\ell_{p, \infty}^{2^m}}\\
		&\leq 2^{1/\beta}. 
	\end{align*}
	Therefore,
	$\|D_{\alpha} (2^{m/p} z)\|_{\ell_{\beta,\infty}^{2^m}} = 2^{m(1/p +\alpha) } \| (i^{-\alpha} z_i)_{i=1}^{2^m} \|_{\ell_{\beta,\infty}^{2^m}}
	\leq 2^{1/\beta} 2^{m/\beta},
	$
	and (\ref{Eq1-Polynomial-decay-Orlicz-F}) is proved.
	
	Next, 
	for any $j \in \{0,1,...,m-5\},$ we write $\psi(m,j) := 2^{m-j-4} (m-j)^{-2}.$   
	Let $W(m)$ be the set given by
	\begin{equation*}
		W(m) = \left\lbrace (x_i)_{i=1}^{2^m} \in \R^{2^m} : x_{\lfloor 2^{j-1} \rfloor +1}^{\ast} \leq 2^{\psi(m,j)} ~\text{for all}~ j \in \{0,1,...,m-5\} \right\rbrace. 
	\end{equation*} 
	Let $x = (x_i)_{i=1}^{2^m} \in B_{\ell_{\beta, \infty}^{2^m}}.$
	We claim that $2^{m/\beta}x \in 4^{-1} 2^{-1/\beta} c(\alpha, p)W(m).$ 
	Fix $j \in \{0,1,...,m-5\}.$
	Then 
	\begin{equation} \label{Eq4-Polynomial-decay-Orlicz-F}
		2^{m/\beta} x_{\lfloor 2^{j-1} \rfloor +1}^{\ast} \leq 2^{m/\beta} (\lfloor 2^{j-1} \rfloor +1)^{-1/\beta}.
	\end{equation}
	If $m=5$ then $j=0.$
	It follows from (\ref{Eq4-Polynomial-decay-Orlicz-F}) that
	\begin{equation} \label{Eq5-Polynomial-decay-Orlicz-F}
		2^{5/\beta} x_1^\ast 
		\leq  	\left( \frac{2^{5/\beta}} {2^{2/25}}\right)  2^{\psi(5,0)}.
	\end{equation}
	Now let us assume that $m \geq 6.$
	If $j=0,$ applying Lemma \ref{Lem-Estimate-using-Calulus} when $t = 2^m,$ we then obtain that 
	\begin{equation} \label{Eq7-Polynomial-decay-Orlicz-F}
		2^{m/\beta} \leq c_1(\beta) 2^{\frac{2^{m-4}}{m^2}} = c_1(\beta) 2^{\psi(m,0)},
	\end{equation}
	where $c_1(\beta) = (22\cdot 16 \log e)^{1/\beta} \left( \frac{1}{\beta} \right) ^{1/\beta} \left( \log \left( \frac{16 \log e}{\beta}\right)  \right)^{2/\beta}.$
	If $1 \leq j \leq m-5,$ by Lemma \ref{Lem-Estimate-using-Calulus} when $t = 2^{m-j},$ we also get that
	\begin{equation} \label{Eq8-Polynomial-decay-Orlicz-F}
		\left(\frac{2^m}{2^{j-1}+1} \right)^{1/\beta} 
		\leq  2^{1/\beta}  2^{(m-j)/\beta}   
		\leq c_2(\beta) 2^{\frac{2^{m-j-4}}{(m-j)^2}} = c_2(\beta) 2^{\psi(m,j)},
	\end{equation}
	where $c_2(\beta) = 2^{1/\beta}c_1(\beta).$
	It follows from (\ref{Eq4-Polynomial-decay-Orlicz-F}), (\ref{Eq7-Polynomial-decay-Orlicz-F}), (\ref{Eq8-Polynomial-decay-Orlicz-F}) that
	$2^{m/\beta} x_{\lfloor 2^{j-1} \rfloor +1}^{\ast} \leq c_2(\beta) 2^{\psi(m,j)}.$
	Hence, the claim was proved.
	The above claim and (\ref{Eq1-Polynomial-decay-Orlicz-F}) imply that 
	\begin{equation} \label{Eq9-Polynomial-decay-Orlicz-F}
		D_{\alpha} \left( \varphi_{M_1}\left( 2^m\right)  B_{\ell_{M_1}^{2^m}} \right) \subseteq  2^{1/\beta} 2^{m/\beta} B_{\ell_{\beta, \infty}^{2^m}} \subseteq 4^{-1}c(\alpha, p) W(m).
	\end{equation}
	By Lemma \ref{Lem-Comnatorial-5}, we can find vectors $y_1,...,y_k \in \ell_{\infty}^{2^m}$ with $k \leq 2^{(2^{m-1}) -1}$ such that
	\begin{equation} \label{Eq10-Polynomial-decay-Orlicz-F}
		W(m) \subseteq \bigcup_{i=1}^{k} \left( y_i + 4 B_{\ell_{\infty}^{2^m}}\right).
	\end{equation}
	We deduce from (\ref{Eq9-Polynomial-decay-Orlicz-F}) and (\ref{Eq10-Polynomial-decay-Orlicz-F}) that 
	\begin{equation} \label{Eq11-Polynomial-decay-Orlicz-F}
		e_{2^{m-1}}\left( D_{\alpha}: \ell_{M_1}^{2^m} \rightarrow \ell_{\infty}^{2^m}\right) 
		\leq \frac{c(\alpha, p)}{\varphi_{M_1}(2^m)}.
	\end{equation}
	Using (\ref{Eq11-Polynomial-decay-Orlicz-F}), we arrive at
	\begin{align} \label{Eq12-Polynomial-decay-Orlicz-F}
		e_{2^{m-1}}\left( D_{\alpha}: \ell_{M_1}^{2^m} \rightarrow \ell_{M_2}^{2^m}\right) 
		&\leq e_{2^{m-1}}\left( D_{\alpha}: \ell_{M_1}^{2^m} \rightarrow \ell_{\infty}^{2^m}\right) \|id: \ell_{\infty}^{2^m} \rightarrow  \ell_{M_2}^{2^m}\| \notag\\
		&\leq c(\alpha, p)  \frac{\varphi_{M_2}(2^m)}{\varphi_{M_1}(2^m)}.
	\end{align}
	The theorem is proved. 
\end{proof}

The next result is an easy consequence of Theorem \ref{Thm-Polynomial-decay-Orlicz-F}.

\begin{cor} \label{Cor1-Polynomial-decay-Orlicz-F}
	Let $n \in \N$ and $\alpha \geq 0.$ 
	Let $D_{\alpha} : \ell_{M_1}^{n} \rightarrow \ell_{M_2}^{n}$ be the linear map defined by 
	\[ D_{\alpha}(x_i)_{i=1}^{n}= \left( \left( \frac{n}{i} \right)^{\alpha} x_i\right) _{i=1}^{n}.\]
	Then,
	\begin{equation} 
		e_{n}\left( D_{\alpha}: \ell_{M_1}^{n} \rightarrow \ell_{M_2}^{n}\right)  \leq c(\alpha, p)  \frac{\varphi_{M_2}(n)}{\varphi_{M_1}(n)},
	\end{equation}
	where $c(\alpha, p) = 4 \cdot 2^{12(1/p +\alpha)} \left(  \log \left( 24(1/p +\alpha)\right) \right)^{2(1/p +\alpha)}   \left( 1/p +\alpha \right) ^{1/p +\alpha}.$ 
\end{cor}
\begin{proof}
	First, we assume that $n \geq 16.$
	Let $m \in \N$ be such that $m \geq 5$ and $2^{m-1} \leq n < 2^m.$
	As $2^m = n\left( n^{-1}2^m\right) ,$ it follows that	 		
	\begin{equation} \label{Eq0-Cor1-Polynomial-decay-Orlicz-F}
		\frac{\varphi_{M_1}( 2^m) }{\varphi_{M_2}(2^m)}
		\geq \left(\frac{2^m}{n} \right)^{-1/p}  \frac{\varphi_{M_1}(n)}{\varphi_{M_2}(n)} 
		\geq 2^{-1/p} \frac{\varphi_{M_1}(n)}{\varphi_{M_2}(n)} .
	\end{equation}
	Theorem \ref{Thm-Polynomial-decay-Orlicz-F} implies that
	\begin{align} \label{Eq3-Cor1-Polynomial-decay-Orlicz-F}
		e_{n} \left( D_{\alpha}: \ell_{M_1}^{n} \rightarrow \ell_{M_2}^{n}\right) 
		&\leq e_{2^{m-1}} \left( D_{\alpha}:\ell_{M_1}^{2^m}\rightarrow \ell_{M_2}^{2^m}\right) \notag\\
		&\leq c_1(p) \frac{\varphi_{M_2}(2^m)}{\varphi_{M_1}(2^m)} \notag\\
		&\leq 2^{1/p} c_1(p) \frac{\varphi_{M_2}(n)}{\varphi_{M_1}(n)},
	\end{align}
	where $c_1(p) = 4 \cdot 2^{11(1/p +\alpha)}   \left(  \log \left( 24(1/p +\alpha)\right) \right)^{2(1/p +\alpha)} \left( 1/p +\alpha \right) ^{1/p +\alpha}.$
	
	Next, let us consider the case $1 \leq n \leq 16.$ 
	Let $(x_i)_{i=1}^n \in B_{\ell_{M_1}^n}.$
	For any fixed $i \in \N,$ we have that $|x_i| \leq M_1^{-1}(1) = 1.$
	It follows that 
	\[ \sum_{i=1}^{n} M_2 \left( \frac{(n/i)^\alpha | x_i|}{ n^\alpha \varphi_{M_2}(n)} \right) 
	\leq \sum_{i=1}^{n} M_2 \left( M_2^{-1}(1/n) \right) 
	=1 . \]		
	Since $\varphi_{M_1} (n) \leq  16^{1/p},$
	we then obtain that
	\begin{equation}\label{Eq1-Cor1-Polynomial-decay-Orlicz-F}
		e_n\left( D_{\alpha}: \ell_{M_1}^{n} \rightarrow \ell_{M_2}^{n}\right)  
		\leq \| D_{\alpha}: \ell_{M_1}^{n} \rightarrow \ell_{M_2}^{n}\|
		\leq  n^\alpha \varphi_{M_2}(n) 
		\leq  16^{1/p+\alpha} \frac{\varphi_{M_2}(n)}{\varphi_{M_1} (n)}.
	\end{equation}
	From (\ref{Eq3-Cor1-Polynomial-decay-Orlicz-F}) and (\ref{Eq1-Cor1-Polynomial-decay-Orlicz-F}), we conclude that
	\begin{equation} \label{Eq2-Cor1-Polynomial-decay-Orlicz-F}
		e_n\left( D_{\alpha}: \ell_{M_1}^{n} \rightarrow \ell_{M_2}^{n}\right)   \leq c_2(\alpha, p)  \frac{\varphi_{M_2}(n)}{\varphi_{M_1} (n)},
	\end{equation}
	where $c_2(\alpha, p) = \max \left\lbrace 2^{1/p} c_1(\alpha, p), 16^{\alpha+1/p} \right\rbrace .$  
\end{proof}


\begin{cor} \label{Cor-Schutt'sThmForOrliczSp-New} 
	Let $n \in \N$ and $id:\ell_{M_1}^{n}\rightarrow \ell_{M_2}^{n}$ the natural embedding.
	Then,
	\begin{equation*}
		c_1(p)  \frac{\varphi_{M_2}(n)}{\varphi_{M_1}(n)} \leq 
		e_{n} \left( id:\ell_{M_1}^{n}\rightarrow \ell_{M_2}^{n}\right) 
		\leq c_2(p) \frac{\varphi_{M_2}(n)}{\varphi_{M_1}(n)},
	\end{equation*}
	where $c_1(p) = 4^{-2} 2^{-12/p }  \left(  \log \left( \frac{24}{p} \right) \right)^{-2/p} \left( \frac{1}{p} \right) ^{-1/p }$ and 
	$c_2(p) = 4 \cdot 2^{12/p }  \left(  \log \left( \frac{24}{p} \right) \right)^{2/p} \left( \frac{1}{p} \right) ^{1/p }.$
\end{cor}
\begin{proof}
	The upper estimate follows from Corollary  \ref{Cor1-Polynomial-decay-Orlicz-F} when $\alpha = 0.$
	To prove the lower estimate, we apply Theorem \ref{Entorpy-of-finite-dim-identity-map-Q} and the upper estimate so that we get
	\begin{align*}
		2^{\frac{1-2n}{n}} &\leq e_{2n}( id:\ell_{M_1}^{n}\rightarrow \ell_{M_1}^{n}) \\
		& \leq e_n( id:\ell_{M_1}^{n}\rightarrow \ell_{M_2}^{n}) e_n( id:\ell_{M_2}^{n}\rightarrow \ell_{M_1}^{n})\\
		& \leq c(p)  \frac{\varphi_{M_1}(n)}{\varphi_{M_2}(n)}  e_n( id:\ell_{M_1}^{n} \rightarrow \ell_{M_2}^{n}),
	\end{align*}
	where $c(p) = 4 \cdot 2^{12/p }  \left(  \log \left( \frac{24}{p} \right) \right)^{2/p} \left( \frac{1}{p} \right) ^{1/p }.$
	The result follows.
\end{proof}

\begin{rem}
	A similar result concerning entropy numbers of embeddings between finite dimensional symmetric Banach spaces can be found in \cite{Sch}.
\end{rem}

Next, let us prove a technical inequality.

\begin{lem} \label{Lem-Inequality-1}
	For any $x \geq 5,$
	\begin{equation}
		\frac{\sqrt{x} \log(2x)}{4 \left( 1+34^{-1} \log(2x) \right) } - \frac{1}{8 \sqrt{x}} > \sqrt{2}.
	\end{equation}
\end{lem}
\begin{proof}
	Since $1+ 34^{-1} \log (2x) \leq \left( \frac{1}{\log 10}+\frac{1}{34} \right)  \log (2x),$ we obtain that
	\begin{equation}  \label{Eq1-Lem-Inequality-1}
		\frac{\sqrt{x} \log(2x)}{4 \left( 1+34^{-1} \log(2x) \right) } - \frac{1}{8 \sqrt{x}}
		\geq \frac{\sqrt{x}}{ 4\left( \frac{1}{\log 10}+\frac{1}{34} \right) }  - \frac{1}{8 \sqrt{x}}.
	\end{equation}
	As the map $t \mapsto f(t) := \frac{\sqrt{t}}{4\left( \frac{1}{\log 10}+\frac{1}{34} \right) }  - \frac{1}{8 \sqrt{t}}$ is increasing for all $t > 0,$ the right-hand side of (\ref{Eq1-Lem-Inequality-1}) is not smaller than $f(5).$
	Since $f(5) > \sqrt{2},$ the result follows.
\end{proof}

Given a set $E$ and $m \in \N$ with $\sharp E \geq m,$ we write
$\mathcal{U}(E,m) = \{U \subseteq E : \sharp U = m \}.$
Recall that $A \Delta B = (A \setminus B) \cup (B\setminus A)$ for any sets $A$ and $B.$

\begin{lem} \label{Lem-Combinatorial1-2}
	Let $k,n \in \N$ be such that $4 \log 2n  \leq  k \leq n/5.$
	Let $s$ be an integer such that
	\begin{equation}
		\frac{k}{\log (2n/k)} < s \leq 1+ \frac{k}{\log (2n/k)}.
	\end{equation}
	Put $E = \{1,2,...,n\}.$ 
	Then there is a set $\mathcal{V}(E,s) \subseteq \mathcal{U}(E,s)$ with the following properties:
	\begin{enumerate}
		\item [\textnormal{(a)}] $\sharp(V_1 \Delta V_2) \geq s$ ~for any distinct $V_1, V_2 \in \mathcal{V}(E,s);$
		\item [\textnormal{(b)}] $\sharp \mathcal{V}(E,s) \geq 2^{k/4}.$
	\end{enumerate}
\end{lem}

\begin{proof}
	Observe that the condition $4 \log 2n  \leq  k \leq n/5$ implies that $ n \geq 168$ and $k \geq 34.$
	We note that $\sharp \mathcal{U}(E,s) = {n \choose s}.$
	Let $V_0 \in \mathcal{U}(E,s).$
	Since $\sharp V_0  + \sharp V_1 = \sharp(V_0 \Delta V_1) + 2 \left[ \sharp(V_0 \cap V_1)\right] ,$  it follows that
	\begin{equation} \label{Eq0-Lem-Combinatorial1-2}
		\{ V_1 \in \mathcal{U}(E,s) : \sharp(V_0 \Delta V_1) \leq s \} = \{ V_1 \in \mathcal{U}(E,s) : \sharp(V_0 \cap V_1) \geq s/2 \}.
	\end{equation}
	Hence, if $s/2 \in \N$ then 
	\begin{equation*} 
		\sharp \{ V_1 \in \mathcal{U}(E,s) : \sharp(V_0 \Delta V_1) \leq s \}  \leq {s \choose s/2} {n-s/2 \choose s/2}.
	\end{equation*}
	Thus there is a set $\mathcal{V}(E,s) \subseteq \mathcal{U}(E,s)$  satisfying the property (a) and that
	\begin{align}  \label{Eq1-Lem-Combinatorial1-2}
		\sharp \mathcal{V}(E,s) 
		&\geq \left. {n \choose s} \right/  \left\lbrace {s \choose s/2}{n-s/2 \choose s/2} \right\rbrace \notag \\
		&= \frac{n(n-1)\cdots (n-s/2+1)}{{s \choose s/2} s(s-1) \cdots (s-s/2+1)} \notag\\
		&\geq \frac{(n-s/2)^{s/2}}{2^s s^{s/2}}.
	\end{align}
	Lemma \ref{Lem-Inequality-1} implies that 
	\begin{equation} \label{Eq2-Lem-Combinatorial1-2}
		\frac{n}{4s} -\frac{1}{8} 
		\geq \frac{1}{4} \cdot\frac{n \log(2n/k)}{k+\log(2n/k)} - \frac{1}{8} 
		\geq \frac{1}{4} \cdot\frac{(n/k) \log(2n/k)}{1+34^{-1}\log(2n/k)} - \frac{1}{8} 
		\geq \left( \frac{2n}{k}\right)^{1/2}.
	\end{equation}
	The estimates (\ref{Eq1-Lem-Combinatorial1-2}) and (\ref{Eq2-Lem-Combinatorial1-2}) give us
	\[ \sharp \mathcal{V}(E,s)  \geq \left( \frac{n}{4s} -\frac{1}{8}  \right)^{s/2}  \geq \left(\frac{2n}{k}\right)^{s/4} \geq 2^{k/4}.\]	
	
	Next, we assume that $s$ is odd.
	Since the following sets are exactly the same:
	$$\{ V_1 \in \mathcal{U}(E,s) : \sharp(V_0 \cap V_1) \geq s/2 \} = \{ V_1 \in \mathcal{U}(E,s) : \sharp(V_0 \cap V_1) \geq (s+1)/2 \},$$
	again by (\ref{Eq0-Lem-Combinatorial1-2}), we have that
	\begin{equation*} 
		\sharp \{ V_1 \in \mathcal{U}(E,s) : \sharp(V_0 \Delta V_1) \leq s \}  \leq {s \choose (s+1)/2} {n-(s+1)/2 \choose (s-1)/2}.
	\end{equation*}
	Thus there is a set $\mathcal{V}(E,s) \subseteq \mathcal{U}(E,s)$  with the property (a) and that
	\begin{align}  \label{Eq3-Lem-Combinatorial1-2}
		\sharp \mathcal{V}(E,s) 
		&\geq \left. {n \choose s} \right/  \left\lbrace {s \choose (s+1)/2} {n-(s+1)/2 \choose (s-1)/2} \right\rbrace \notag \\
		&\geq \frac{n(n-1)\cdots (n-(s+1)/2+1)}{{s \choose (s+1)/2} s(s-1) \cdots (s-(s+1)/2+1)} \notag\\
		&\geq \frac{(n-s/2)^{(s+1)/2}}{2^{s+1} s^{(s+1)/2}}.
	\end{align}
	Using (\ref{Eq2-Lem-Combinatorial1-2}) and (\ref{Eq3-Lem-Combinatorial1-2}), we then obtain that
	\[\sharp \mathcal{V}(E,s)  \geq \left( \frac{n}{4s} -\frac{1}{8}  \right)^{(s+1)/2}  \geq \left(\frac{2n}{k}\right)^{(s+1)/4} \geq 2^{k/4}. \]
	The lemma is proved.
\end{proof}

The following theorem is an extension of Sch\"{u}tt's result (see \cite{Sch}, Theorem 1) when the natural embedding is defined on Orlicz sequence spaces.  
We note that the corresponding upper estimate will be included later (see Theorem \ref{Thm-Schutt'sThmForOrliczSp-New2}) because it is a consequence of our main result (Theorem \ref{GeneralEdm&Yu-Orlicz-sp-Q}) in the next section .
In the special case when $M_1(t) = t^p$ and  $M_2(t) = t^q$ for any $0<p<q<\infty,$ some estimates of constants can also be found in \cite{Gued&Litv}. 

\begin{thm} \label{Thm-Schutt'sThmForOrliczSp-New1} 
	Let $k,n \in \N.$
	Then,
	\begin{equation*}
		e_k(id:\ell_{M_1}^n\rightarrow \ell_{M_2}^n) 
		\geq  
		\begin{cases}
			2^{-1/p}	&~~\text{if}~~ 1\leq k <  \log 2n  ,\\
			c(p)  \displaystyle\frac{\varphi_{M_2}  \left(k/\log (2n/k) \right)   }{\varphi_{M_1}\left( k/\log (2n/k) \right)  } 				&~~\text{if}~~ \log 2n  \leq k \leq n,
		\end{cases}
	\end{equation*}
	where $c(p) =4^{-2}  2^{-19/p }  \left(  \log \left( \frac{24}{p} \right) \right)^{-2/p} \left( \frac{1}{p} \right) ^{-1/p }.$
	
\end{thm}
\begin{proof}
	The case $1 \leq k <  \log 2n $ is elementary.
	We assume first that $ \log 2n  \leq k \leq n/20.$
	Set $m = 4k.$
	Then, $4\log 2n  \leq m \leq n/5.$ 
	Let $E = \{1,2...,n\}$ and let $s$ be a positive integer such that 
	\[ \frac{m}{\log\left(2n/m \right) } < s \leq 1+ \frac{m}{\log\left( 2n/m \right) }.\] 
	By Lemma \ref{Lem-Combinatorial1-2},  there is a set $\mathcal{V}(E,s) := \{V \subseteq E : \sharp V = s \}$ with the following properties:
	\begin{enumerate}
		\item [\textnormal{(a)}] $\sharp(V_1 \Delta V_2) \geq s$ ~for any distinct $V_1, V_2 \in \mathcal{V}(E,s);$
		\item [\textnormal{(b)}] $\sharp \mathcal{V}(E,s) \geq 2^{m/4}.$
	\end{enumerate}
	It is easy to verify that for every $V \in \mathcal{V}(E,s),$ 
	$~\frac{1}{\varphi_{M_1}(s) }\chi_{V} \in B_{\ell_{M_1}^n},$
	where $\chi_{V}$ is the characteristic function of $V.$
	Furthermore, by the property (a), for any distinct $V_1,V_2 \in \mathcal{V}(E,s),$ 
	\begin{equation} \label{Eq1-Lem-Schutt'sThmForOrliczSp-New1} 
		\left\|\frac{1}{\varphi_{M_1}(s) }\chi_{V_1} -\frac{1}{\varphi_{M_1}(s) }\chi_{V_2}  \right\|_{\ell_{M_2}^n} 
		\geq \frac{\varphi_{M_2}(s)}{\varphi_{M_1}(s) }.
	\end{equation}
	We note that the condition $ \log 2n  \leq k \leq n/20$ implies that $n \geq 168$ and $k \geq 9.$
	As \[ \frac{k}{\log\left(2n/k \right) } \leq s < 1+ \frac{4k}{\log\left(n/2k \right) } \leq  \frac{10k}{\log\left(2n/k \right) },\]
	the following inequalities hold:
	\[ \varphi_{M_1}(s)\leq 10^{1/p} \varphi_{M_1}\left(  \frac{k}{\log (2n/k)} \right) ~~\text{and}~~ 
	\varphi_{M_2}\left(  \frac{k}{\log (2n/k)} \right) \leq \varphi_{M_2}(s).\]
	Since $\sharp \left\lbrace \frac{1}{\varphi_{M_1}(s) }\chi_{V}: V \in \mathcal{V}(E,s)  \right\rbrace \geq 2^{m/4} > 2^{k -1},$ 
	it follows from (\ref{Eq1-Lem-Schutt'sThmForOrliczSp-New1}) and the previous remarks that
	\begin{align*}
		e_{k} (id:\ell_{M_1}^n\rightarrow \ell_{M_2}^n) 
		&\geq  2^{1-1/p} f_{k} (id:\ell_{M_1}^n\rightarrow \ell_{M_2}^n) \\
		&\geq  2^{1-1/p} 2^{-1} \frac{\varphi_{M_2}(s)}{\varphi_{M_1}(s) }\\
		&\geq   20^{-1/p} \frac{\varphi_{M_2}\left( k/\log (2n/k) \right) }{\varphi_{M_1}\left( k/\log (2n/k)  \right) }.
	\end{align*}
	Now we assume that $n/20 \leq k \leq n.$ 
	Then, $n \log (2n/k) \leq  (20 \log 40)k.$ 
	Applying Corollary \ref{Cor-Schutt'sThmForOrliczSp-New}, 
	we obtain that
	\begin{align*}
		e_k(id:\ell_{M_1}^n \rightarrow \ell_{M_2}^n) 
		&\geq e_n(id:\ell_{M_1}^n \rightarrow \ell_{M_2}^n) \\
		&\geq c_1(p)  \frac{\varphi_{M_2}(n)}{\varphi_{M_1}(n)}\\
		&\geq  c_1(p)  (20 \log 40)^{-1/p} \frac{\varphi_{M_2} \left( k/\log (2n/k) \right)  }{\varphi_{M_1}\left( k/\log (2n/k) \right)  },
	\end{align*}
	where $c_1(p) = 4^{-2}  2^{-12/p }  \left(  \log \left( \frac{24}{p} \right) \right)^{-2/p} \left( \frac{1}{p} \right) ^{-1/p }.$
\end{proof}


The following useful estimate was proved in \cite{E&N2010}:
\begin{lem} \label{Lem-Edmunds-Netrusov-1}
	Let $m \in \N$ and let $\{E_i\}_{i=0}^m$ be a family of non-empty disjoint sets such that $\sharp E_i \leq 2^{m+2^i}.$
	Denote by $\mathcal{F}$ the family of all sequences $(F_i)_{i=0} ^m$ such that $F_i \subseteq E_i$ and $\sum_{i=0}^{m} \frac{\sharp F_i}{2^{m-i}} \leq 1.$
	Then, 
	\begin{equation}
		\log \left( \sharp \mathcal{F} \right) \leq 2^{m+3}.
	\end{equation}
\end{lem}

Using the idea of Edmunds and Netrusov in \cite{E&N2010} in a general setting, we obtain the estimates of a particular diagonal map acting between Orlicz sequence spaces. 

\begin{lem} \label{Lemma-Before-main-Orlicz-seq-spaces-Q}
	Let $m \in \N$ be such that $m \geq 5$ and let $E_0,E_1,....,E_{m+1}$ be non-empty disjoint sets such that $\sharp E_i \leq 2^{m+2^i}$ for all $i = 0,1,...,m;$ it is possible that $\sharp E_{m+1} = \infty.$
	Assume that $\varphi_{M_1}/\varphi_{M_2}$ is a non-decreasing function.
	Define a sequence $(\omega_i)_{i=0}^{m+1}$ by
	\[
	\omega_i =
	\begin{cases}
		\displaystyle \frac{\varphi_{M_1}(2^{m-i})}{\varphi_{M_2}(2^{m-i})} &~~\text{if}~~ i = 0,1,...,m,\\
		1  &~~\text{if}~~ i = m+1.
	\end{cases}
	\]
	Put $E = \bigcup_{i=0}^{m+1} E_i$  
	and let $T_{\omega}: \ell_{M_1}(E) \rightarrow \ell_{M_2}(E)$ be the linear map defined  by 
	$T_{\omega}(x) = \left( \omega_i  x_i \right)_{i=0}^{m+1}$
	for  $x = (x_i)_{i=0}^{m+1} \in \ell_{M_1}(E)$ and $x_i \in \ell_{M_1}(E_i).$
	Then, 
	\begin{equation}
		e_{10 \cdot 2^m}\left( T_{\omega}: \ell_{M_1}(E) \rightarrow \ell_{M_2}(E)\right)  \leq c(p) ,
	\end{equation}
	where  $c(p) = 4 \cdot 2^{12/p }  \left(  \log \left( \frac{24}{p} \right) \right)^{2/p} \left( \frac{1}{p} \right) ^{1/p } .$
\end{lem}

\begin{proof}
	For any family of sets $\{F_i \}_{i=0}^{m}$ such that $F_i \subseteq E_i$ for all $i =0,1,...,m$ and $ \sum_{i=0}^{m} (\sharp F_i)2^i  \leq 2^{m},$ 
	we put $F = \cup_{i=0}^{m}F_i $ and let $id_F : \ell_{M_1}(F) \rightarrow \ell_{M_2}(F)$ be the natural embedding.
	Set $\vartheta := e_{2^m}(id_F : \ell_{M_1}(F) \rightarrow \ell_{M_2}(F)).$
	Denote by $\Gamma_{F}$ an optimal $\vartheta$-net, i.e., a net of cardinality $2^{2^m-1}$ of the set $B_{\ell_{M_1}(F)}$ in $\ell_{M_2}(F).$
	Let $\Gamma = \bigcup \Gamma_{F},$ where the union is taken over the family of all sets $\{F_i \}_{i=0}^m$ as above. 
	By Lemma \ref{Lem-Edmunds-Netrusov-1}, 
	\begin{equation} \label{Eq1-Lemma-Before-main-Orlicz-seq-spaces-Q}
		\sharp \Gamma \leq 2^{2^{m+3}} 2^{2^m-1} = 2^{9\cdot 2^m-1}.
	\end{equation}
	
	We will prove that $\Gamma$ is an $\eta$-net of $T_{\omega}(B_{\ell_{M_1}(E)})$ in $\ell_{M_2}(E),$
	where 
	\begin{equation} \label{Eq2-Lemma-Before-main-Orlicz-seq-spaces-Q}
		\eta^{p} =    \omega_0^p e_{ 2^m}^{p}\left( id: \ell_{M_1}^{2^m} \rightarrow \ell_{M_2}^{2^m}\right)  + 1.
	\end{equation}
	Let $x = (x_i)_{i=0}^{m+1} \in B_{\ell_{M_1}(E)},$ where, for each $i \in \{0,1,...,m+1\}, ~x_i = (\xi_{i,j})_{j\in E_i}.$ 
	For any $i \in \{0,1,...,m+1\},$ we define $F_i \subseteq E_i$ by
	\[  F_i = \left\lbrace j \in E_i :  M_1(|\xi_{i,j}|)  \geq 2^{i-m} \right\rbrace,\]	
	in addition, we set $y_i = x_i \chi_{F_i},$ where $\chi$ denotes the characteristic function, and put $z_i = x_i - y_i.$
	It is clear that $F_{m+1} = \emptyset.$ 
	Moreover, since $\|x\|_{\ell_{M_1}(E)} \leq 1,$ it follows that
	\[ 1 \geq 	\sum_{i=0}^{m} \sum_{j \in F_i} M_1(|\xi_{i,j}|)
	\geq \sum_{i=0}^{m} (\sharp F_i)2^{i-m},\]
	and hence $\sharp F \leq 2^m.$		
	To obtain the desired result in (\ref{Eq2-Lemma-Before-main-Orlicz-seq-spaces-Q}), we will show that $\|T_{\omega}(z_i)_{i=0}^{m+1} \|_{ \ell_{M_2}(E)} \leq 1.$		
	Let us fix $i \in \{0,1,...,m+1\}$ and $j \in E_i \setminus F_i.$
	We first consider  the case that $i = m+1.$ 
	Since $ \varphi_{M_1}/\varphi_{M_2}$ is non-decreasing, 
	$|\xi_{m+1,j}| = M_1^{-1} \left(  M_1 (|\xi_{m+1,j}|)\right) \leq M_2^{-1} \left(  M_1 (|\xi_{m+1,j}|)\right).$
	This implies that
	\begin{equation}  \label{Eq3-Lemma-Before-main-Orlicz-seq-spaces-Q}
		M_2 \left(\omega_{m+1}|\xi_{m+1,j}| \right) \leq M_1 (|\xi_{m+1,j}|).
	\end{equation}
	Next, we assume that $i \in \{0,1,...,m\}$.
	As $j \notin F_i,$ we have that $|\xi_{i,j}| \leq M_{1}^{-1}(2^{i-m}) = 1/\varphi_{M_1}(2^{m-i}).$ 
	If $|\xi_{i,j}| \neq 0$ then there is $t \geq 2^{m-i}$ such that $|\xi_{i,j}| = 1/\varphi_{M_1}(t) = M_1^{-1}(1/t).$
	It follows that 
	\begin{equation}  \label{Eq4-Lemma-Before-main-Orlicz-seq-spaces-Q}
		M_2 \left( \omega_i |\xi_{i,j}| \right) 
		\leq M_2 \left( \frac{\varphi_{M_1}(t)}{\varphi_{M_2}(t)} \frac{1}{ \varphi_{M_1}(t)} \right)
		= \frac{1}{t}
		= M_1(|\xi_{i,j}|).
	\end{equation}
	Note that the inequality (\ref{Eq4-Lemma-Before-main-Orlicz-seq-spaces-Q}) is also true when $|\xi_{i,j}|  = 0.$
	We deduce from (\ref{Eq3-Lemma-Before-main-Orlicz-seq-spaces-Q}) and (\ref{Eq4-Lemma-Before-main-Orlicz-seq-spaces-Q}) that
	\[ \sum_{i=0}^{m+1} \sum_{j \in E_i\setminus F_i} M_2 \left( \omega_i |\xi_{i,j}| \right) \leq 1.\]
	Therefore, $\|T_{\omega}(z_i)_{i=0}^{m+1} \|_{ \ell_{M_2}(E)} \leq 1.$			
	It follows from (\ref{Eq2-Lemma-Before-main-Orlicz-seq-spaces-Q}) and Theorem \ref{Thm-Polynomial-decay-Orlicz-F} (when $\alpha = 0$)  that
	\begin{align*} 
		e_{10\cdot 2^m}(T_{\omega}: \ell_{M_1}(E) \rightarrow \ell_{M_2}(E)) 
		&\leq \left(  \omega_0^p \left( c(p) \frac{\varphi_{M_2}(2^m)}{\varphi_{M_1}(2^m)}\right)^p  +1 \right)^{1/p}\\
		&= (c(p)^p +1)^{1/p},
	\end{align*}		
	where  $c(p) =4 \cdot 2^{11/p }  \left(  \log \left( \frac{24}{p} \right) \right)^{2/p} \left( \frac{1}{p} \right) ^{1/p }.$
	The proof is completed.
\end{proof}	

\begin{cor}  \label{Cor-Construct-E}
	Let $m \in \N$ be such that $m \geq 5.$
	Assume that $\varphi_{M_1}/\varphi_{M_2}$ is a non-decreasing function.
	Let $(\alpha_i)_{i=1}^\infty$ be a non-increasing sequence of non-negative real numbers such that
	\begin{equation} \label{Eq1-Cor-Construct-E}
		\alpha_1  \leq \frac{\varphi_{M_1}(2^m)}{\varphi_{M_2}(2^m)} ~~~\text{~~~and~~~}~~~ \alpha_{2^{m+2^{i-1}}}  \leq \frac{\varphi_{M_1}(2^{m-i})}{\varphi_{M_2}(2^{m-i})} ~~\text{~~for all}~i \in \{1,2,...,m \}.
	\end{equation}
	Let $T_\alpha \in \mathcal{B}( \ell_{M_1},\ell_{M_2})$ be the diagonal operator generated by $(\alpha_i)_{i=1}^\infty.$ 
	Then,
	\[e_{10 \cdot 2^m}(T_\alpha: \ell_{M_1} \rightarrow \ell_{M_2}) \leq c(p), 
	\]
	where  $c(p) = 4 \cdot 2^{12/p }  \left(  \log \left( \frac{24}{p} \right) \right)^{2/p} \left( \frac{1}{p} \right) ^{1/p } .$
\end{cor}
\begin{proof}
	Construct the sets $E_0, E_1,...,E_{m+1}$ as follows:
	\begin{equation*}
		\begin{cases}
			E_0 &= \left\lbrace 1,2,...,2^{m+1}-1 \right\rbrace, \\
			E_i &= \left\lbrace 2^{m+2^{i-1}},...,2^{m+2^{i}}-1 \right\rbrace ~\text{for}~ i=1,2,...,m,\\
			E_{m+1} &= \left\lbrace 2^{m+2^{m}}, 2^{m+2^{m}}+1,... \right\rbrace.
		\end{cases}
	\end{equation*}
	Put $E = \bigcup_{i=0}^{m+1} E_i.$
	The desired upper estimate follows from Lemma \ref{Lemma-Before-main-Orlicz-seq-spaces-Q} and the fact that $e_{10 \cdot 2^m}(T_\alpha) \leq e_{10 \cdot 2^m}\left( T_{\omega}\right),$
	where $T_{\omega}: \ell_{M_1}(E) \rightarrow \ell_{M_2}(E)$ is the diagonal map defined in Lemma \ref{Lemma-Before-main-Orlicz-seq-spaces-Q}.
\end{proof}
\section{Main results and some consequences}

Now we are ready to prove our main result.

\begin{thm}\label{GeneralEdm&Yu-Orlicz-sp-Q}
	Let $(\alpha_i)_{i=1}^\infty$ be a non-increasing sequence of non-negative real numbers with $\alpha_1 = \alpha_{k}$ for some $k \in \N.$
	Assume that $\varphi_{M_1}/\varphi_{M_2}$ is a non-decreasing function.
	Set 
	\begin{equation}  \label{Eq0-General-Edm&Yu-In-Orlicz-Q}
		\Theta(k) := \max_{k \leq n \leq 2^{k-1}} \alpha_n 
		\frac{\varphi_{M_2} \left(  k/\log (2n/k)   \right) }{\varphi_{M_1}\left( k/\log (2n/k)  \right)} .
	\end{equation}
	Let $T_\alpha \in \mathcal{B}( \ell_{M_1},\ell_{M_2})$ be the diagonal operator generated by $(\alpha_i)_{i=1}^\infty.$ 
	Then, 
	\[  c_1(p)  \Theta(k) \leq e_k(T_\alpha: \ell_{M_1} \rightarrow \ell_{M_2}) \leq c_2(p)  \Theta(k),\]
	where $c_1(p) = 4^{-2}  2^{-19/p }  \left(  \log \left( \frac{24}{p} \right) \right)^{-2/p} \left( \frac{1}{p} \right) ^{-1/p }$ and 
	$c_2(p) = 4 \cdot 2^{19/p }  \left(  \log \left( \frac{24}{p} \right) \right)^{2/p} \left( \frac{1}{p} \right) ^{1/p }.$
\end{thm}

\begin{proof}
	To prove the estimate from below, let us fix $n \in \N$ with $k \leq n \leq 2^{k-1}.$ 
	Let $P_n: \ell_{M_1} \rightarrow \ell_{M_1}$ be the projection onto the first $n$ coordinates.
	Then, by Theorem \ref{Thm-Schutt'sThmForOrliczSp-New1} , the following estimates hold:
	\begin{equation*} 
		e_k(T_\alpha) \geq e_k(T_{\alpha}\circ P_n) 
		\geq \alpha_n e_k(id: \ell_{M_1}^n \rightarrow \ell_{M_2}^n)  
		\geq c_1(p) \alpha_n 	
		\frac{\varphi_{M_2} \left(  k/\log (2n/k)  \right) }{\varphi_{M_1}\left(  k/\log (2n/k)  \right)} ,
	\end{equation*}
	where $c_1(p) = 4^{-2}  2^{-19/p }  \left(  \log \left( \frac{24}{p} \right) \right)^{-2/p} \left( \frac{1}{p} \right) ^{-1/p }.$
	The desired lower estimate follows.
	
	Next, let us prove the upper estimate.  
	Let $(\sigma_j)_{j=1}^\infty$ be a sequence defined by 
	\[ 
	\sigma_j = 
	\begin{cases}
		\displaystyle 2^{-7/p}\frac{\varphi_{M_1}(k)}{\varphi_{M_2}(k)}    &\text{if }  1 \leq j \leq k,\\
		\displaystyle  2^{-7/p}\frac{\varphi_{M_1} \left(  k/\log (2j/k) \right) }{\varphi_{M_2}\left(  k/\log (2j/k)  \right)}     &\text{if } k \leq j \leq 2^{k-1}, \\
		\displaystyle  2^{-7/p}\frac{\varphi_{M_1} \left(  k/\log (2^k/k)  \right) }{\varphi_{M_2}\left(  k/\log (2^k/k)  \right)}     &\text{if }  j \geq 2^{k-1}.
	\end{cases}
	\]
	Since $\varphi_{M_1}/\varphi_{M_2}$ is non-decreasing, $(\sigma_j)_{j=1}^\infty$ is a non-increasing sequence.		 					
	Define $T_\sigma: \ell_{M_1} \rightarrow \ell_{M_2}$ by $$T_\sigma ((x_j)_{j=1}^{\infty}) = (\sigma_j x_j)_{j=1}^{\infty}.$$
	Suppose first that $k \geq 2^9.$ 
	Let $m \in \N$ be such that $m \geq 5$ and $ 2^{m+4} \leq k < 2^{m+5}.$ 
	We will show that $(\sigma_j)_{j=1}^\infty$ satisfies (\ref{Eq1-Cor-Construct-E}).
	As $k \leq 2^{m+4},$ we have that
	\begin{equation}
		\sigma_1 = 2^{-7/p}\frac{\varphi_{M_1}(k)}{\varphi_{M_2}(k)}
		\leq 2^{-2/p} \frac{\varphi_{M_1}(2^m)}{\varphi_{M_2}(2^m)} \leq \frac{\varphi_{M_1}(2^m)}{\varphi_{M_2}(2^m)}.
	\end{equation}
	Next, let us fix $i \in \{1,2,...,m\}.$
	If $i \in \{1,2\}$ then $2^{m+2^i-1} \leq 2^{m+4} \leq k \leq  2^{m+7-i}.$
	This implies that
	\begin{equation}
		\sigma_{2^{m+2^i-1}} = 2^{-7/p}\frac{\varphi_{M_1}(k)}{\varphi_{M_2}(k)} \leq \frac{\varphi_{M_1}(2^{m-i})}{\varphi_{M_2}(2^{m-i})}.
	\end{equation}
	Now we consider the case that $i \in \{3,4,...,m\}.$
	It is easy to check that $k \leq 2^{m+2^i-1} \leq 2^{k-1}.$
	Since $\log(2^{m+2^i}/k) \geq \log(2^{2^i - 5}) \geq 2^{i-2}$ and $2^{m-i} \log(2^{m+2^i}/k) \leq 2^m \leq k,$ it follows that
	\begin{equation}
		\sigma_{2^{m+2^i-1}} = 2^{-7/p}\frac{\varphi_{M_1} \left(  k/\log (2^{m+2^i}/k)  \right) }{\varphi_{M_2}\left( k/\log (2^{m+2^i}/k)  \right)} \leq 
		2^{-7/p} \frac{\varphi_{M_1}(2^{m+5}/ 2^{i-2})}{\varphi_{M_2}(2^{m-i})} \leq
		\frac{\varphi_{M_1}(2^{m-i})}{\varphi_{M_2}(2^{m-i})}.
	\end{equation}
	Therefore $(\sigma_j)_{j=1}^\infty$ satisfies (\ref{Eq1-Cor-Construct-E}).
	Corollary \ref{Cor-Construct-E} yields
	\begin{equation}
		e_k(T_\sigma: \ell_{M_1} \rightarrow \ell_{M_2}) \leq e_{10\cdot 2^m}\left( T_\sigma: \ell_{M_1} \rightarrow \ell_{M_2}\right)  \leq c_{2}(p) ,
	\end{equation}
	where $c_2(p) = 4 \cdot 2^{12/p }  \left(  \log \left( \frac{24}{p} \right) \right)^{2/p} \left( \frac{1}{p} \right) ^{1/p } .$
	Since the following estimate holds:  
	\[ e_k(T_{\alpha}: \ell_{M_1} \rightarrow \ell_{M_2})  \leq e_k\left( 2^{7/p} \Theta(k)T_\sigma: \ell_{M_1} \rightarrow \ell_{M_2} \right) ,\]
	we obtain that $e_k(T_{\alpha}: \ell_{M_1} \rightarrow \ell_{M_2})   \leq 2^{7/p}c_2(p) \Theta(k).$
	
	Finally, we assume that $1 \leq k \leq 2^9.$ 
	Let us first estimate the norm of $T_{\alpha}.$
	Fix $(x_i)_{i=1}^\infty \in B_{\ell_{M_1}}.$
	Then $|x_i| \leq M_1^{-1}(1) = 1/\varphi_{M_1}(1)$ for all $i \in \N.$
	If $|x_i| \neq 0$ then there exists $t \geq 1$ such that $|x_i| = 1/\varphi_{M_1}(t) = M_1^{-1}(1/t).$
	It follows that 
	\begin{equation} \label{Eq6-General-Edm&Yu-In-Orlicz-Q}
		M_2 \left( \frac{|\alpha_i x_i|}{\alpha_{k}} \right) 
		\leq M_2 \left( \frac{\varphi_{M_1}(t)}{\varphi_{M_2}(t)} \frac{1}{ \varphi_{M_1}(t)} \right) 
		= \frac{1}{t}
		= M_1(|x_i|) . 
	\end{equation}
	We note that the inequality (\ref{Eq6-General-Edm&Yu-In-Orlicz-Q}) is also true when $|x_i| = 0.$
	Thus, we have
	$\sum_{i=1}^{\infty} M_2 \left( \frac{|\alpha_i x_i|}{\alpha_k} \right) \leq 1.$
	Hence,
	$e_k(T_{\alpha}: \ell_{M_1} \rightarrow \ell_{M_2})   \leq \|T_{\alpha}\| \leq \alpha_{k}.$ 
	The desired upper estimate follows from the following inequalities: 
	\[
	\Theta(k) 
	\geq \alpha_{k}  \frac{\varphi_{M_2}(k)}{\varphi_{M_1}(k)}  
	\geq    2^{-9/p} e_k(T_{\alpha}: \ell_{M_1} \rightarrow \ell_{M_2})  .
	\]
	The proof is completed.
\end{proof}
It might be of interest to state the following extension of Sch\"{u}tt's result (see \cite{Sch}, Theorem 1) when the natural embedding is defined on Orlicz sequence spaces.  

\begin{thm} \label{Thm-Schutt'sThmForOrliczSp-New2} 
	Let $k,n \in \N.$ If $\varphi_{M_1}/\varphi_{M_2}$ is a non-decreasing function, then there are positive constants $c_1(p)$ and  $c_2(p)$ such that
	\begin{equation*}
		c_1(p) A(n,k) \leq  e_k(id:\ell_{M_1}^n\rightarrow \ell_{M_2}^n) \leq c_2(p) A(n,k) 
	\end{equation*}
	where 	
	\[
	\begin{cases}
		1	&~~\text{if}~~ 1\leq k \leq  \log 2n  ,\\
		\displaystyle\frac{\varphi_{M_2}  \left(k/\log (2n/k) \right)   }{\varphi_{M_1}\left( k/\log (2n/k) \right)  } 				&~~\text{if}~~ \log 2n  \leq k \leq n,\\
		2^{-k/n} \dfrac{\varphi_{M_2} (n)}{\varphi_{M_1} (n)}	&~~\text{if}~~  k \geq n .
	\end{cases}
	\]
\end{thm}
\begin{proof}
	If $k \geq n $, then the desired result is a special case of Theorem 4.2 in \cite{Kaewtem1}. 
	For the case $1\leq k \leq  \log 2n$, the upper estimate follows from the inequality $ e_k(id:\ell_{M_1}^n\rightarrow \ell_{M_2}^n) \leq \left\| id:\ell_{M_1}^n\rightarrow \ell_{M_2}^n \right\|.$
	So, in view of Theorem \ref{Thm-Schutt'sThmForOrliczSp-New1}, it remains to give  the upper estimate when $\log 2n  \leq k \leq n.$
	Let $T_\alpha \in \mathcal{B}( \ell_{M_1},\ell_{M_2})$ be the diagonal operator generated by $(\alpha_i)_{i=1}^\infty$ where
	\[
	\alpha_1 = \alpha_2 = ... = \alpha_k = ... = \alpha_n =1 ~~\text{~~and~~}~~ \alpha_i = 0 ~\text{~for all~}~ i \geq n+1.
	\]
	Theorem \ref{GeneralEdm&Yu-Orlicz-sp-Q} implies that
	\begin{align*}
		e_k(id:\ell_{M_1}^n\rightarrow \ell_{M_2}^n) 
		&\leq e_k(T_\alpha:\ell_{M_1}\rightarrow \ell_{M_2}) \\
		&\leq c(p) \max_{k \leq m \leq n} \displaystyle\frac{\varphi_{M_2}  \left(k/\log (2m/k) \right)   }{\varphi_{M_1}\left( k/\log (2m/k) \right)  }\\
		&= c(p)  \displaystyle\frac{\varphi_{M_2}  \left(k/\log (2n/k) \right)   }{\varphi_{M_1}\left( k/\log (2n/k) \right)  }
	\end{align*}
	
\end{proof}

Next, for each $m \in \N,$ we write
\begin{equation} \label{Eq-Lambda-notation}
	\Lambda(m) := \max_{s \in \{1,...,m\}} \alpha_{m2^{s-1}} \frac{\varphi_{M_2} \left( m/s \right) }{\varphi_{M_1}\left(  m/s  \right)}. 
\end{equation}

\begin{prop} \label{Rem-Equiv-of-Theta-Q}
	Let $(\alpha_i)_{i=1}^\infty$ be a non-increasing sequence of non-negative real numbers. 
	Then, for each $k \in \N,$
	\begin{equation*}
		2^{-2/p}\Lambda(k) \leq \Theta(k) \leq 2^{1/p}\Lambda(k),
	\end{equation*}
	where $\Theta(k) $ was defined in Theorem  \ref{GeneralEdm&Yu-Orlicz-sp-Q} .
\end{prop}
\begin{proof}
	First, let us prove the lower estimate.
	Fix $s \in \{1,...,k\}.$ 
	If $k \leq k2^{s-1} \leq 2^{k-1},$ we immediately obtain that $\Theta(k) \geq \alpha_{k2^{s-1}} \frac{\varphi_{M_2} \left( k/s \right) }{\varphi_{M_1}\left(  k/s  \right)}. $
	Now we consider the case that $k2^{s-1} \geq 2^{k-1}.$
	This implies that $\varphi_{M_2}  \left( k/\log (2^k/k)  \right) \geq \varphi_{M_2} \left( k/s \right).$
	Since $\log (2^k/k) \geq k/4 \geq s/4,$ we also have $\varphi_{M_1}  \left( k/\log (2^k/k)  \right) \leq  2^{2/p}\varphi_{M_1} \left( k/s \right).$
	We deduce that  
	\[ \Theta(k) \geq \alpha_{2^{k-1}} \frac{\varphi_{M_2} \left( k/\log (2^k/k) \right) }{\varphi_{M_1}\left( k/\log (2^k/k)  \right)}
	\geq 2^{-2/p}\alpha_{k2^{s-1}} \frac{\varphi_{M_2} \left( k/s \right) }{\varphi_{M_1}\left(  k/s  \right)}.  \]
	The desired lower estimate follows.
	
	Next, to prove the upper estimate, let $n \in \N$ be such that $k \leq n \leq 2^{k-1}.$
	Set $m = \lfloor \log(2n/k) \rfloor.$
	Then, $1 \leq m \leq k.$
	As $2^{-1} \log (2n/k) \leq m\leq \log(2n/k)$ and $k2^{m-1} \leq n,$ we thus obtain that
	\[ \Lambda(k) \geq \alpha_{k2^{m-1}} \frac{\varphi_{M_2} \left( k/m \right) }{\varphi_{M_1}\left(  k/m \right)}
	\geq 2^{-1/p}\alpha_n   \frac{\varphi_{M_2} \left(  k/\log (2n/k)  \right) }{\varphi_{M_1}\left( k/\log (2n/k) \right)}. \]
\end{proof}
Next, we apply the main result to obtain entropy estimates for diagonal operators generated by sequences satisfying the doubling condition.

\begin{cor} \label{Cor-Doubling-Orlicz-sp-Q}
	Let $(\alpha_i)_{i=1}^\infty$ be a non-increasing sequence of non-negative real numbers. 
	Assume that $(\alpha_i)_{i=1}^\infty$ satisfies the condition that there is a constant $C \geq 1$
	such that for any $i \in \N$
	\begin{equation}\label{Doubling-condition-Q}
		\alpha_i \leq C \alpha_{2i}
	\end{equation}
	Let $T_\alpha \in \mathcal{B}( \ell_{M_1},\ell_{M_2})$ be the diagonal operator generated by $(\alpha_i)_{i=1}^\infty.$ 
	If  $ \varphi_{M_1}/\varphi_{M_2}$ is a non-decreasing function, then,  for any $k \in\N,$ 
	\[ c_1(p) \Lambda(k) \leq e_k(T_\alpha) \leq c_2(p, C) \Lambda(k),\] 
	where  $c_1(p)$ and $c_2(p, C)$ are constants given by
	$c_1(p) = 4^{-2}  2^{-21/p }  \left(  \log \left( \frac{24}{p} \right) \right)^{-2/p} \left( \frac{1}{p} \right) ^{-1/p }$
	and $c_2(p, C) = 4 C  2^{22/p } 2^{14 C} \left(  \log \left( 24(1/p +C)\right) \right)^{2(1/p +C)}   \left( 1/p +C \right) ^{1/p +C} .$
\end{cor}

\begin{proof}
	First, let us prove the upper estimate.
	We suppose that $k \geq 2.$
	Choose $m \in \N$ such that $2m \leq k < 2(m+1).$
	Let $ (\sigma_i)_{i=1}^{\infty}$ and $(\mu_i)_{i=1}^{\infty}$ be sequences defined by
	\[ 
	\sigma_i = 
	\begin{cases}
		\alpha_{m}  	&~\text{if}~~  i \leq m,\\
		\alpha_i 	&~\text{if}~~ i> m,
	\end{cases}
	~~~\text{~~~and~~~}~~~
	\mu_i = 
	\begin{cases}
		\displaystyle \left( \frac{4m}{i}\right)^{C}  \alpha_{4m} 	&~\text{if}~~  i \leq m,\\
		0 	&~\text{if}~~ i> m.
	\end{cases}
	\]
	Define $T_{\sigma}: \ell_{M_1} \rightarrow \ell_{M_2} $ and $T_{\mu} : \ell_{M_1} \rightarrow \ell_{M_2}$ by
	$T_{\sigma} (x_{i})_{i=1}^\infty= (\sigma_i x_i)_{i=1}^\infty$ and $T_{\mu} (x_{i})_{i=1}^\infty = (\mu_i x_i)_{i=1}^\infty,$ respectively.
	By Theorem \ref{GeneralEdm&Yu-Orlicz-sp-Q} and Proposition \ref{Rem-Equiv-of-Theta-Q}, we have that
	\begin{equation} \label{Eq2-Diagonal2-op-Orlicz-Q}
		c_1(p)  \Lambda(m) \leq e_m(T_{\sigma}: \ell_{M_1} \rightarrow \ell_{M_2}) \leq c_2(p)  \Lambda(m),
	\end{equation}
	where 
	$c_1(p) = 4^{-2}  2^{-21/p }  \left(  \log \left( \frac{24}{p} \right) \right)^{-2/p} \left( \frac{1}{p} \right) ^{-1/p }$ and 
	$c_2(p) = 4 \cdot 2^{20/p }  \left(  \log \left( \frac{24}{p} \right) \right)^{2/p} \left( \frac{1}{p} \right) ^{1/p }.$
	We note that the condition (\ref{Doubling-condition-Q}) implies that 
	\begin{equation}
		i^C\alpha_i  \leq C j^C \alpha_j ~\text{~~~for~~}~~1 \leq i \leq j < \infty.
	\end{equation}
	Hence, for any $i \leq m,$
	\begin{equation}  \label{Eq3-Diagonal2-op-Orlicz-Q}
		\alpha_i \leq C \left( \frac{4m}{i}\right)^{C} \alpha_{4m}   = C\mu_i.
	\end{equation}	
	It can be shown that
	$\Lambda (m) \leq C 4^{1/p+C}  \Lambda (k)$ and $\varphi_{M_2}(m)/\varphi_{M_1}(m) \leq 4^{1/p} \varphi_{M_2}(k)/\varphi_{M_2}(k).$
	Corollary \ref{Cor1-Polynomial-decay-Orlicz-F}, (\ref{Eq2-Diagonal2-op-Orlicz-Q}) and (\ref{Eq3-Diagonal2-op-Orlicz-Q}) give us 
	\begin{align} \label{Eq4-Diagonal2-op-Orlicz-Q}
		e_{2m}^p(T_{\alpha}: \ell_{M_1} \rightarrow \ell_{M_2}) 
		&\leq   e_m^p(T_{\sigma} :\ell_{M_1} \rightarrow \ell_{M_2}) + C^{p} e_{m}^{p}(T_{\mu} :\ell_{M_1} \rightarrow \ell_{M_2})  \notag \\
		&\leq  \left[  c_2(p)  \Lambda(m)\right] ^p + C^{p}\left[ 4^C c_3(p,C)  \alpha_{4m} \frac{\varphi_{M_2}(m)}{\varphi_{M_1}(m)} \right]^p \notag\\
		&\leq  \left[  c_2(p)  C 4^{1/p+C}  \Lambda (k)\right] ^p + C^{p}\left[ 4^{1/p+C}c_3(p,C)  \alpha_{k} \frac{\varphi_{M_2}(k)}{\varphi_{M_1}(k)} \right]^p \notag\\
		&\leq  \left[  c_2(p)  C 4^{1/p+C}  \Lambda (k)\right] ^p + C^{p}\left[ 4^{1/p+C}c_3(p,C) \Lambda(k) \right]^p \notag\\
		&\leq \left( c_2^{p}(p) + c_3^{p} (p,C) \right)   \left( C 4^{1/p+C} \Lambda(k)\right) ^p,
	\end{align}
	where $c_3(p, C) = 4 \cdot 2^{12(1/p +C)} \left(  \log \left( 24(1/p +C)\right) \right)^{2(1/p +C)}   \left( 1/p +C\right) ^{1/p +C}.$ 
	It follows from (\ref{Eq4-Diagonal2-op-Orlicz-Q}) that
	\begin{equation}
		e_{k}(T_{\alpha}: \ell_{M_1} \rightarrow \ell_{M_2}) 
		\leq e_{2m}(T_{\alpha}: \ell_{M_1} \rightarrow \ell_{M_2}) 
		\leq c_4(p,  C)  \Lambda(k),
	\end{equation}
	where $\displaystyle c_4(p, C) = C 4^{1/p +C}   \left( c_2^{p}(p) + c_3^{p}(p,C)  \right)^{1/p}.$
	As $e_1(T_{\alpha}: \ell_{M_1} \rightarrow \ell_{M_2}) \leq \|T_{\alpha}\| \leq \alpha_1  =  \Lambda(1),$ the proof of the upper estimate is completed.
	
	The lower estimate follows from (\ref{Eq2-Diagonal2-op-Orlicz-Q}) (when $m=k$) and the inequality $e_k(T_{\alpha}) \geq e_k(T_{\sigma}).$
\end{proof}	



\begin{prop}\label{Ex1-Polyn-decay-Orlicz-Q}
	Let $k \in \N$ and $\theta>0.$ 
	Let $T_\alpha : \ell_{M_1} \rightarrow \ell_{M_2}$ be the linear map defined by 
	\[  T_\alpha(x_i)_{i=1}^\infty =  \left( \alpha_i x_i\right) _{i=1}^{\infty}, ~~\text{where}~~ 
	\alpha_i = i^{-\theta}
	~~\text{for all}~~ i \in \N.
	\]
	If  $\varphi_{M_1}/\varphi_{M_2}$  is a non-decreasing function,  then there are positive constants $c_1(p)$ and $c_2(p, \theta)$ such that
	\[  c_1(p) \alpha_k \frac{\varphi_{M_2}(k)}{\varphi_{M_1}(k)} \leq e_k(T_\alpha: \ell_{M_1} \rightarrow \ell_{M_2}) \leq c_2(p, \theta) \alpha_k   \frac{\varphi_{M_2}(k)}{\varphi_{M_1}(k)}.
	\]
\end{prop}
\begin{proof}
	%
	The result follows from Corollary \ref{Cor-Doubling-Orlicz-sp-Q} and the fact that
	$x^{1/p} \leq c(p, \theta)
	2^{\theta x}$ for all $x \geq 1.$
\end{proof}

\begin{rem}
	In the context of Lebesgue sequence spaces, the above result was due to Carl in 1981 (see \cite{Carl1981-diagonal}, Proposition 1 and Proposition 2).
\end{rem}


\begin{prop} \label{Example-exp-log}
	Let $k \in \N, \beta >0$ and $0<\vartheta<1.$ 
	For each $i \in \N,$ let 
	\[ \alpha_i =  \exp \left( -\beta \log^{\vartheta} (i+1)\right) .\]
	Let $T_{\alpha} \in \mathcal{B}( \ell_{M_1}, \ell_{M_2})$ be the diagonal operator generated by $(\alpha_i)_{i=1}^\infty.$ 
	If  $\varphi_{M_1}/\varphi_{M_2}$ is a non-decreasing function, then there are positive constants $c_1$ and $c_2,$ depending only on $p, \beta$ and $\vartheta,$ such that 
	\begin{equation*}
		c_1 \alpha_k \frac{\varphi_{M_2} \left( k/ \log^{1-\vartheta} (k+1)\right)}{\varphi_{M_1} \left( k/ \log^{1-\vartheta} (k+1)\right)} \leq e_{k}(T_\alpha)  \leq c_2 \alpha_k \frac{\varphi_{M_2} \left( k/ \log^{1-\vartheta} (k+1)\right)}{\varphi_{M_1} \left( k/ \log^{1-\vartheta} (k+1)\right)}.
	\end{equation*}
\end{prop}
\begin{proof}
	It is clear that $(\alpha_i)_{i=1}^\infty$ is a non-increasing sequence.
	Let  $1 \leq i \leq j < \infty.$
	By the Mean Value Theorem, there is a real number $t_0$ such that $\log (i+1)<t_0< \log(j+1)$ and
	$
	\log^{\vartheta} (j+1) - \log^{\vartheta} (i+1)
	= \vartheta t_0^{\vartheta -1} \left( \log(j+1) - \log (i+1) \right). 
	$
	This implies that
	\begin{align*}
		\frac{\alpha_i}{\alpha_j} &= \exp \left(\beta \left[ \log^{\vartheta}(j+1) -  \log^{\vartheta}(i+1) \right]  \right) \\
		&\leq \exp\left( \beta  \vartheta  \log \left( \frac{j+1}{i+1}\right)   \right) \\
		&\leq \left(\frac{j}{i} \right)^{\beta \vartheta \log e} .
	\end{align*}
	Thus $(\alpha_i)_{i=1}^\infty$ satisfies the condition (\ref{Doubling-condition-Q}) with $C=2^{\beta \vartheta \log e}.$
	Now let us prove the lower estimate.
	Put $m = \lfloor \log^{1-\vartheta} (k+1) \rfloor.$
	Then, $1 \leq m \leq k.$
	By the Mean Value Theorem, we have that 
	\begin{align*}
		\frac{\alpha_k}{\alpha_{k2^{m-1}}} 
		&\leq \exp \left( \beta \left[ \left( m+\log (k+1)\right)^{\vartheta}  - \log^{\vartheta}(k+1) \right]  \right) \\
		&\leq \exp \left( \beta \vartheta  m /\log^{1-\vartheta}(k+1) \right) \\
		&\leq \exp \left( \beta \vartheta  \right). 
	\end{align*}
	It follows from Corollary \ref{Cor-Doubling-Orlicz-sp-Q} that
	$$e_{k}(T_\alpha:\ell_{M_1} \rightarrow \ell_{M_2}) 
	\geq c_1 \alpha_{k2^{m-1}} \frac{\varphi_{M_2}(k/m)}{\varphi_{M_1}(k/m)}
	\geq c_2\alpha_k \frac{\varphi_{M_2} \left( k/ \log^{1-\vartheta} (k+1)\right)}{\varphi_{M_1} \left( k/ \log^{1-\vartheta} (k+1)\right)}.
	$$
	
	Next, to prove the upper estimate, let us fix $s \in \{1,...,k\}.$ 
	By Corollary \ref{Cor-Doubling-Orlicz-sp-Q}, the proof will be completed if we show that there exists a constant $C_0 :=C_0(p, \beta, \vartheta)>0$ such that 
	\vspace{-0.3cm}
	\begin{equation} \label{Eq2-Example-exp-log}
		\alpha_{k2^{s-1}} \frac{\varphi_{M_2}(k/s)}{\varphi_{M_1}(k/s)} \leq C_0 \alpha_k \frac{\varphi_{M_2} \left( k/ \log^{1-\vartheta} (k+1)\right)}{\varphi_{M_1} \left( k/ \log^{1-\vartheta} (k+1)\right)}.
	\end{equation} 
	Set $C_1 := C_1 (p,\beta,\vartheta) = \left( 1+2^{2-\vartheta} (p \beta \vartheta)^{-1} \right) ^{1/\vartheta}.$		
	First, we assume that $1 \leq s \leq C_1\log^{1-\vartheta}(k+1).$
	Then,
	$$\frac{\varphi_{M_1} \left( k/ \log^{1-\vartheta} (k+1)\right)}{\varphi_{M_2} \left( k/ \log^{1-\vartheta} (k+1)\right)}
	\leq C_1^{1/p} \frac{\varphi_{M_1} \left( k/ C_1\log^{1-\vartheta} (k+1)\right)}{\varphi_{M_2} \left( k/ C_1\log^{1-\vartheta} (k+1)\right)}
	\leq C_1^{1/p}  \frac{\varphi_{M_1} \left( k/ s\right)}{\varphi_{M_2} \left( k/ s\right)}.$$
	Since $(\alpha_i)_{i=1}^\infty$ is non-increasing, we arrive at
	\begin{equation} \label{Eq3-Example-exp-log}
		\alpha_{k2^{s-1}} \frac{\varphi_{M_2}(k/s)}{\varphi_{M_1}(k/s)} 
		\leq C_1^{1/p}  \alpha_k \frac{\varphi_{M_2} \left( k /\log^{1-\vartheta} (k+1)\right)}{\varphi_{M_1} \left( k/ \log^{1-\vartheta} (k+1)\right)}.
	\end{equation}
	Now we consider the case that $ C_1\log^{1-\vartheta}(k+1) \leq s \leq k.$
	As the following estimates hold: $\varphi_{M_2}(k/s) \leq  \varphi_{M_2} \left( k/ \log^{1-\vartheta} (k+1)\right)$ and $\varphi_{M_1} \left( k/ \log^{1-\vartheta} (k+1)\right) \leq \left( s/ \log^{1-\vartheta} (k+1)\right)^{1/p} \varphi_{M_1}(k/s),$
	to obtain (\ref{Eq2-Example-exp-log}), it is enough to show that
	\begin{equation} \label{Eq9-Example-exp-log}
		\alpha_{k2^{s-1}} \left(\frac{s}{\log^{1-\vartheta}(k+1)} \right)^{1/p} \leq C_0\alpha_k.
	\end{equation}
	Let $f$ be the function defined by 
	\[ f(t) = t^{1/p} \exp\left( -\beta \log^{\vartheta} (k2^{t-1}+1)\right)  ~~\text{for all}~~ t \geq C_1\log^{1-\vartheta}(k+1).
	\]
	Then,
	$$f^{\prime}(t) =  t^{1/p} \exp\left( -\beta \log^{\vartheta} (k2^{t-1}+1)\right) \left[ \frac{1}{p t} - \frac{\beta \vartheta}{\log^{1-\vartheta}(k2^{t-1}+1) }\left( \frac{k2^{t-1}}{k2^{t-1}+1}\right)  \right] .$$
	Using the facts that $0 < \vartheta < 1$ and $t \geq C_1\log^{1-\vartheta}(k+1)\geq C_1,$ we obtain that
	\[  t/(t + \log(k+1))^{1-\vartheta} 
	\geq 2^{\vartheta-1} \min \left\lbrace  t^{\vartheta}, t/\log^{1-\vartheta}(k+1) \right\rbrace 
	\geq 2/(p \beta \vartheta).
	\]
	The above inequalities yield
	\[ \frac{t}{\log^{1-\vartheta}(k2^{t-1}+1) }\left( \frac{k2^{t-1}}{k2^{t-1}+1}\right) 
	\geq
	\frac{t}{(t + \log(k+1))^{1-\vartheta}}\cdot\frac{1}{2}
	\geq \frac{1}{p \beta \vartheta}.
	\]
	This implies that $f$ is decreasing for all $t \geq C_1\log^{1-\vartheta}(k+1).$
	It follows that
	\begin{align*} \label{Eq8-Example-exp-log}
		\hspace{-0.75cm}
		s^{1/p} \exp\left( -\beta \log^{\vartheta} (k2^{s-1}+1)\right) 
		&\leq C_1^{1/p} \log^{(1-\vartheta)/p}(k+1) \exp\left( -\beta \log^{\vartheta} \left( k 2^{C_1 \log^{1-\vartheta}(k+1)-1} +1\right)  \right) \\
		&\leq C_1^{1/p} \log^{(1-\vartheta)/p}(k+1) \exp\left( -\beta \log^{\vartheta} \left( k +1\right)  \right).
	\end{align*}
	The above inequalities imply (\ref{Eq9-Example-exp-log}) with $C_0 = C_1^{1/p}.$
\end{proof}

\begin{rem}
	Let $0<p < q <\infty.$ 
	If $M_1(t) = t^p$ and $M_2(t) = t^q,$  the diagonal operator in Example  \ref{Example-exp-log} was studied by K\"{u}hn (see \cite{Ku2}, Theorem 4.3).
	The estimate of $e_{k}(T:\ell_{p} \rightarrow \ell_{q})$ from below was sharp; however, the upper estimate was not optimal.
	The sharp two-sided estimates were given by Edmunds and Netrusov (see \cite{E&N2010}, Example 14). 
\end{rem}


\begin{prop}\label{Ex1-LOG-decay-Orlicz-Q}
	Let $k \in \N$ and $\theta >0.$
	Let $T_\alpha : \ell_{M_1} \rightarrow \ell_{M_2}$ be the linear map defined by 
	$T_\alpha(x_i)_{i=1}^\infty =  \left( \alpha_i x_i\right) _{i=1}^{\infty},$ where 
	$\alpha_i = \left( \log(i+1)\right)^{-\theta} $
	for all $i \in \N.$
	If  $\varphi_{M_1}/\varphi_{M_2}$ is a non-decreasing function, then there are positive constants $c_1(p, \theta)$ and $c_2(p, \theta)$ such that
	\[ c_1(p, \theta) \Phi(k) \leq e_k(T_\alpha: \ell_{M_1} \rightarrow \ell_{M_2}) \leq c_2(p, \theta)  \Phi(k),\]
	where
	$\Phi(k) := \displaystyle \max \left\lbrace \alpha_{2^{s-1}} \frac{\varphi_{M_2}(k/s)}{\varphi_{M_1}(k/s)} :  s \in \N ~\text{and}~ \log 2k \leq s \leq k \right\rbrace.$
\end{prop}

\begin{proof}
	Since the map $t \mapsto t^{-1} \log(t+1)$ is decreasing for $t \geq 1,$  for any $1\leq i \leq j < \infty,$ we have that
	\[\frac{\alpha_i}{\alpha_j} = \left( \frac{\log(j+1)}{\log (i+1)} \right)^{\theta} \leq \left( \frac{j}{i}\right)^{\theta} . \]
	Thus $(\alpha_i)_{i=1}^\infty$ satisfies the condition (\ref{Doubling-condition-Q}) with $C= 2^\theta.$
	First, let us prove the estimate from below.
	For $\log 2k \leq s \leq k,$ one can see that 
	$\log\left( k2^{s-1}+1 \right) \leq \log k + \log \left(2^{s-1}+1 \right) \leq 2 \log \left(2^{s-1}+1 \right) ,$
	and hence 
	\[ \alpha_{k2^{s-1}}= \left( \frac{1}{\log\left( k2^{s-1}+1 \right)}\right)^{\theta} 
	\geq \left( \frac{1}{2\log \left(2^{s-1}+1 \right) } \right)^{\theta}
	=2^{-\theta} \alpha_{2^{s-1}}. 
	\]
	This implies that
	$\Lambda(k) \geq 2^{-\theta} \Phi(k).$
	The lower estimate follows from Corollary \ref{Cor-Doubling-Orlicz-sp-Q}.
	
	To prove the upper estimate, let us fix $s \in \{1,2,...,k\}.$
	If $\log 2k \leq s \leq k,$ then 
	\begin{equation}
		\alpha_{k2^{s-1}} \frac{\varphi_{M_2}(k/s)}{\varphi_{M_1}(k/s)} 
		\leq \alpha_{2^{s-1}} \frac{\varphi_{M_2}(k/s)}{\varphi_{M_1}(k/s)} 
		\leq \Phi(k).
	\end{equation}
	For the case $1\leq s \leq \log 2k,$ we also get that
	\begin{equation}
		\alpha_{k2^{s-1}} \frac{\varphi_{M_2}(k/s)}{\varphi_{M_1}(k/s)} 
		\leq  \alpha_{k} \frac{\varphi_{M_2}(k/\log 2k)}{\varphi_{M_1}(k/\log 2k)} 
		\leq   \Phi(k).
	\end{equation}
	By Corollary \ref{Cor-Doubling-Orlicz-sp-Q}, the upper inequality is proved.
\end{proof}

\begin{rem}
	Let $\theta>0, 0<p<q<\infty, M_1(t) = t^p$ and $M_2(t) = t^q.$ 
	Let $T_\alpha:\ell_p \rightarrow \ell_q$ be the diagonal operator defined in Example \ref{Ex1-LOG-decay-Orlicz-Q}. 
	Then there are positive constants $c_1(\theta,p,q)$ and $c_2(\theta,p,q)$ such that
	\[ 
	c_1(\theta,p,q) \phi(k)  \leq e_k(T_\alpha:\ell_p \rightarrow \ell_q) \leq c_2(\theta,p,q) \phi(k),
	\]
	where
	\[ 
	\phi(k):=
	\begin{cases}
		k^{-\theta} &\text{if}~ \theta \leq 1/p-1/q, \\
		k^{1/q-1/p} (\log 2k)^{1/p-1/q-\theta} &\text{if}~ \theta \geq 1/p-1/q.
	\end{cases}
	\]
	If $\theta = 1/p-1/q$ then the above estimates were proved in  \cite{Cobos&Kuhn&Schonbek} and  \cite{Ku2}  for the case that $0 < p < q=\infty$ and $0<p<q \leq \infty,$ respectively. 
	For general $\theta >0,$ two-sided estimates were given in \cite{Ku3} where the restriction $1 \leq p < q \leq \infty$ is required. 
\end{rem}

\section{Final remarks}

\begin{defn}
	Let $K\geq 1.$
	A function $\psi:[0,\infty) \rightarrow [0,\infty)$ is {\it $K$-increasing} if $\psi(s) \leq K\psi(t)$ for any $0\leq s\leq t.$
\end{defn}

The following result is known for experts in the area (see, for instance, the proof of Lemma 1 on page 76 of \cite{E&N1998}).

\begin{thm} \label{Prop-Edmund-Netrusov-exist-concave}
	Let $\psi$ be a non-negative 
	function defined on $[0,\infty)$ such that $\psi(0) = 0$ and $ \psi(t)>0$ for $t>0.$ 
	If $\psi$ and the map $t \mapsto \frac{t}{\psi(t)}$ are $K$-increasing for some $K \geq 1,$ then there are a continuous, increasing and concave function $\Psi$ defined on $[0,\infty),$ and positive constants $c_1,c_2$ which are independent of $\psi$ and $K$ such that, for all $t>0,$
	\[c_1\psi(t) \leq \Psi(t) \leq c_2 K^2 \psi(t). \]
\end{thm}

\begin{rem} \label{Rmk-l_M-is-Quasi-B-Sp}
	Let $M$ be a non-negative  function defined on $[0,1]$ such that $M(0) = 0, M(1) = 1$ and $M(t) >0$ for all $t>0.$
	If $M$ the map $t \mapsto \frac{M(t)}{t^p}$ are increasing for some $p  \in (0,1],$ then there are an Orlicz function $N$ on $[0,1]$ and positive constants $c_1(p)$ and $c_2(p),$ depending only on $p,$  such that  
	\[  c_1(p)M(t) \leq N(t^p) \leq c_2(p) M(t).\]
	Therefore, it is possible to give a generalisation of our main result by relaxing the convexity of $M(t^{1/p})$.
\end{rem}

\begin{rem}
	~
	\begin{enumerate}
		\item Let $0<p\leq 1$ and $0<q \leq \infty.$ 
		In \cite{Gued&Litv}, Gu\'{e}don and Litvak showed that there are absolute positive constants $c_1$ and $c_2$ such that, for $\log n \leq k \leq n,$
		\begin{equation} \label{Eq1-Final-rmk}
			c_1^{1/p}  \zeta(n,k) 
			\leq e_k(id:\ell_p^n \rightarrow \ell_q^n) 
			\leq c_2^{1/p} \left( \frac{1}{p} \log \left(  \frac{2}{p}\right)  \right)^{1/p-1/q}  \zeta(n,k) ,
		\end{equation}
		where $\zeta(n,k) = \left( k^{-1} \log (2n/k) \right)^{1/p-1/q}.$
		Let us compare the constants in (\ref{Eq1-Final-rmk}) when $q \rightarrow \infty$ with our result. 
		We can see that the main different term in the upper estimate is only $\left( \log(2/p) \right)^{1/p}$ but our constant in the lower estimate is quite far from the optimal one. 
		However, we consider the more general spaces.
		\item A key result used to prove the lower estimate of the main theorem (Theorem \ref{GeneralEdm&Yu-Orlicz-sp-Q}) is Sch\"{u}tt's theorem in a general setting.
		In view of Theorem 4.2 in \cite{E&N1998}, it is possible to obtain the lower estimate of Theorem \ref{GeneralEdm&Yu-Orlicz-sp-Q} in the context of symmetric $p$-Banach spaces.
		In contrast,  the upper estimate of Theorem \ref{GeneralEdm&Yu-Orlicz-sp-Q} cannot be obtained easily even in the case of diagonal operators acting between Lorentz sequence spaces, for example, from $\ell_{p,r}$ to $\ell_{p,s}$ where $0<p<r<s \leq \infty.$
		A reason for this difficulty is that entropy numbers do not behave well under real interpolations (see \cite{E&N2010}).
	\end{enumerate}
\end{rem}








\section*{Acknowledgements}

	We would like to acknowledge the School of Mathematics, University of Bristol for support of this work.
Furthermore, we thank Professor David Edmunds for careful reading of the manuscript and helpful suggestions.
Finally, we thank the referees for providing valuable comments which helped to improve the content of this paper.
The first-named 
author was supported by Mahidol Wittayanusorn School and the Ministry of Science and Technology of Thailand.


\vspace{0.25cm}

(Thanatkrit Kaewtem) {\sc Department of Mathematics and Computing Science, Mahidol Wittayanusorn School, Nakhon Pathom 73170, Thailand}

Email address: {\tt thanatkrit.ktm@mwit.ac.th}
\vspace{0.25cm}

(Yuri Netrusov) {\sc School of Mathematics, Faculty of Science, University of Bristol, Bristol BS8 1TW, UK }

Email address: {\tt y.netrusov@bristol.ac.uk}


\begin{thebibliography}{[1]}
	\bibitem{Abra&Steg}
	M. Abramowitz and I.A. Stegun, {\it Handbook of Mathematical Functions}, Dover, New York, 1972.
	
	\bibitem{Aoki}
	T. Aoki, {\it Locally bounded linear topological spaces}, Proc. Imp. Acad. \textbf{18} (1942), no. 10, 588--594 .
	
	\bibitem{Birman&Solomyak1}
	M.S. Birman and M.Z. Solomyak, {\it Piecewise polynomial approximations of functions of the class ${W}_p^\alpha$}, Mat. Sb. (N.S.), \textbf{73(115)}:3 (1967), 331--355; Math. USSR-Sb., \textbf{2}:3 (1967), 295--317.
	
	
	
	\bibitem{Carl1981-diagonal}
	B. Carl, {\it Entropy numbers of diagonal operators with application to eigenvalue problems}, J. Approx. Theory \textbf{32} (1981), 135--150.
	
	\bibitem{C&S}
	B. Carl and I. Stephani, {\it Entropy, compactness and the approximation of operators}, Cambridge University Press, New York, 1990.
	
	\bibitem{Cobos&Kuhn&Schonbek}
	F. Cobos and T. K\"{u}hn and T. Schonbek, {\it Compact embeddings of {B}r\'{e}zis-{W}ainger type}, Rev. Mat. Iberoam. \textbf{22}(1) (2006), 305--322.
	
	\bibitem{E&N1998}
	D.E. Edmunds and Yu. Netrusov, {\it Entropy numbers of embedding of {S}obolev spaces in {Z}ygmund spaces}, Studia Math. \textbf{128} (1998), 71--102.
	
	\bibitem{E&N2010}
	D.E. Edmunds and Yu. Netrusov, {\it Entropy numbers and interpolation}, Math. Ann. \textbf{351} (2010), 963--977.
	
	\bibitem{E&T}
	D.E. Edmunds and H. Triebel, {\it Function Spaces, entropy numbers, differential operators}, Cambridge University Press, Cambridge, 1996.
	
	\bibitem{Gor&Ko&Sch}
	Y. Gordon and H. K\"{o}nig and C. Sch\"{u}tt, {\it Geometric and probabilistic estimates for entropy and approximation numbers of operators}, J. Approx. Theory \textbf{49} (1987), 219--239.
	
	\bibitem{Gued&Litv}
	O. Gu\'{e}don and A.E. Litvak, {\it Euclidean projections of a $p$-convex body}; in: Geometric Aspects of Functional Analysis, Lecture Notes in Mathematics vol. 1745, Springer-Verlag, Berlin, 2000, pp. 95--108.
	
	\bibitem{Haroske&Triebel1994}
	D. Haroske and H. Triebel, {\it Entropy numbers in weighted function spaces and eigenvalue distribution of some degenerate pseudodifferential operators I}, Math. Nachr. \textbf{167} (1994), 131--156.
	
	\bibitem{Kaewtem1}
	T. Kaewtem, {\it Entropy numbers in $\gamma$-Banach spaces}, Math. Nachr. \textbf{290} (2017), 2879--2889.
	
	\bibitem{Kol&Ti}
	A.N. Kolmogorov and V.M. Tichomirov, {\it $\varepsilon$-entropy and $\varepsilon$-capacity of sets in function spaces}, Uspekhi Mat. Nauk \textbf{14}(2) (1959), 3--86 (in Russian); English transl.: Amer. Math. Soc. Transl. Ser. 2, \textbf{17} (1961),  277--364.
	
	\bibitem{Ku3}
	T. K\"{u}hn, {\it Entropy numbers of diagonal operators of logarithmic type}, Georgian Math. J. \textbf{8}(2) (2001), 307--318.
	
	\bibitem{Ku2}
	T. K\"{u}hn, {\it Entropy numbers of general diagonal operators}, Rev. Mat. Complut. \textbf{18} (2005), 479--491.
	
	\bibitem{Ku&Leo&Sick&Skr1}
	T. K\"{u}hn, H.-G. Leopold, W. Sickel, and L. Skrzypczak, {\it Entropy numbers of embeddings of weighted Besov spaces}, Constr. Approx. \textbf{23} (2006), 61--77.
	
	\bibitem{Ku&Leo&Sick&Skr2}
	T. K\"{u}hn, H.-G. Leopold, W. Sickel, and L. Skrzypczak, {\it Entropy numbers of embeddings of weighted Besov spaces II}, Proc. Edinb. Math. Soc. (2) \textbf{49} (2006), 331--359.
	
	\bibitem{Ku&Leo&Sick&Skr3}
	T. K\"{u}hn, H.-G. Leopold, W. Sickel, and L. Skrzypczak, {\it Entropy numbers of embeddings of weighted Besov spaces III. Weights of logarithmic type}, Math. Z. \textbf{255} (2007), no. 1, 1--15.
	
	\bibitem{Lind&Tzaf}
	J. Lindenstrauss and L. Tzafriri, {\it Classical {B}anach spaces I ({S}equence spaces)}, Springer-Verlag, Berlin, 1977.
	
	\bibitem{Orlicz}
	W. Orlicz, {\it \"{U}ber eine gewisse {K}lasse von {R}\"{a}umen vom {T}ypus {B}}, Bull. Intern. Acad. Pol. \textbf{8} (1932), 207--220.
	
	\bibitem{P0}
	A. Pietsch, {\it Operator ideals}, Deutscher Verlag der Wissenschaften, Berlin, 1978.
	
	\bibitem{P1}
	A. Pietsch, {\it Eigenvalues and s-numbers}, Cambridge University Press, Leipzig, 1987.
	
	\bibitem{Rolewicz}
	S. Rolewicz, {\it On a certain class of linear metric spaces}, Bull. Acad. Polon. Sci. \textbf{5} (1957), 471--473.
	
	\bibitem{Sch}
	C. Sch\"{u}tt, {\it Entropy numbers of diagonal operators between symmetric {B}anach spaces}, J. Approx. Theory \textbf{40} (1984), 121--128.
	
	
\end{thebibliography}
\end{document}